\documentclass[oneside,a4paper]{amsart}

\usepackage[T1]{fontenc}
\usepackage[utf8]{inputenc}
\usepackage{amsmath, amsthm, amsfonts, amssymb}
\usepackage{mathtools}
\usepackage{stmaryrd}

\usepackage{url}
\usepackage[all]{xy}
\usepackage{placeins} 
\usepackage[usenames, dvipsnames]{xcolor} 
\usepackage{euscript} 
\usepackage{bbm} 
\usepackage{enumitem} 
\usepackage{comment}
\usepackage{colonequals}
\usepackage{nicefrac}
\usepackage{float}
\usepackage{CJKutf8}

\usepackage{tabularray}

\usepackage{hyperref}
\hypersetup{
    colorlinks, linkcolor=NavyBlue,
    citecolor=OliveGreen, urlcolor=RedOrange
}
\usepackage[nameinlink,capitalise,noabbrev]{cleveref} 
\usepackage{todonotes}
\setuptodonotes{fancyline, color=green!30}
\usepackage{faktor}

\usepackage{tikz-cd}
\usepackage{tikz}
\usetikzlibrary{arrows,backgrounds,
    decorations.pathreplacing,
    decorations.pathmorphing
}
\usetikzlibrary{matrix,arrows}

\makeatletter
\DeclareFontFamily{OMX}{MnSymbolE}{}
\DeclareSymbolFont{MnLargeSymbols}{OMX}{MnSymbolE}{m}{n}
\SetSymbolFont{MnLargeSymbols}{bold}{OMX}{MnSymbolE}{b}{n}
\DeclareFontShape{OMX}{MnSymbolE}{m}{n}{
    <-6>  MnSymbolE5
   <6-7>  MnSymbolE6
   <7-8>  MnSymbolE7
   <8-9>  MnSymbolE8
   <9-10> MnSymbolE9
  <10-12> MnSymbolE10
  <12->   MnSymbolE12
}{}
\DeclareFontShape{OMX}{MnSymbolE}{b}{n}{
    <-6>  MnSymbolE-Bold5
   <6-7>  MnSymbolE-Bold6
   <7-8>  MnSymbolE-Bold7
   <8-9>  MnSymbolE-Bold8
   <9-10> MnSymbolE-Bold9
  <10-12> MnSymbolE-Bold10
  <12->   MnSymbolE-Bold12
}{}

\let\llangle\@undefined
\let\rrangle\@undefined
\DeclareMathDelimiter{\llangle}{\mathopen}%
                     {MnLargeSymbols}{'164}{MnLargeSymbols}{'164}
\DeclareMathDelimiter{\rrangle}{\mathclose}%
                     {MnLargeSymbols}{'171}{MnLargeSymbols}{'171}
\makeatother

\newcommand{\euscr}[1]{\EuScript{#1}} 
\newcommand{\acat}{\euscr{A}} 
\newcommand{\ccat}{\euscr{C}} 
\newcommand{\Hom}{\textnormal{Hom}} 
\newcommand{\Ext}{\textnormal{Ext}} 
\newcommand{\map}{\textnormal{map}} 
\newcommand{\abeliangroups}{\euscr{A}b} 
\newcommand{\ZZ}{\mathbf{Z}}
\newcommand{\NN}{\mathbf{N}}

\newcommand{\thesphere}{\mathbf{S}}

\newcommand{\spectra}{\euscr{S}p}

\newcommand{\MU}{\mathrm{MU}} 

\newcommand{\Hrm}{\mathrm{H}} 
\newcommand{\Mod}{\euscr{M}\mathrm{od}}

\newcommand{\fieldp}{\mathbf{F}_{p}}

\theoremstyle{plain}

\newtheorem{theorem}{Theorem}[section]

\newtheorem{lemma}[theorem]{Lemma}

\newtheorem{proposition}[theorem]{Proposition}

\newtheorem{corollary}[theorem]{Corollary}
\newtheorem{conjecture}[theorem]{Conjecture}

\newtheorem*{theorem*}{Theorem}

\theoremstyle{definition}
\newtheorem{example}[theorem]{Example}

\newtheorem{definition}[theorem]{Definition}

\newtheorem{remark}[theorem]{Remark}
\newtheorem{recollection}[theorem]{Recollection}

\newtheorem{notation}[theorem]{Notation}

\newtheorem*{remark*}{Remark}
\newtheorem*{interpretation*}{Interpretation}
\newtheorem*{definition*}{Definition}
\newtheorem*{conjecture*}{Conjecture}
\newtheorem*{notation*}{Notation}

\numberwithin{equation}{section}

\makeatletter
  \def\subsection{\@startsection{subsection}{1}%
  \z@{.7\linespacing\@plus\linespacing}{.5\linespacing}%
  {\normalfont\bfseries\centering}}
\makeatother

\let\oldtocsection=\tocsection
\let\oldtocsubsection=\tocsubsection
\let\oldtocsubsubsection=\tocsubsubsection
\renewcommand{\tocsection}[2]{\hspace{0em}\oldtocsection{#1}{#2}}
\renewcommand{\tocsubsection}[2]{\hspace{1em}\oldtocsubsection{#1}{#2}}
\renewcommand{\tocsubsubsection}[2]{\hspace{2em}\oldtocsubsubsection{#1}{#2}}

\calclayout

\newcommand{\mfrak}{\mathfrak{m}} 

\newcommand{\BP}{\mathrm{BP}}
\newcommand{\K}{\mathrm{K}}
\newcommand{\BPn}[1]{\mathrm{BP} \langle #1 \rangle}

\DeclareMathOperator{\Ker}{Ker}
\DeclareMathOperator{\cof}{cofib}

\begin{document}
\title[The monochromatic Hahn-Wilson conjecture]{The monochromatic Hahn-Wilson conjecture}

\author{David Jongwon Lee}
\address{Department of Mathematics, Northwestern University, Evanston, IL, USA}
\email{davidlee@northwestern.edu}

\author{Piotr Pstr\k{a}gowski}
\address{Hakubi Center and Research Institute for Mathematical Sciences, Kyoto University, Japan}
\email{pstragowski.piotr@gmail.com}

\begin{abstract}
We prove the $K(n)$-local analogue of the Hahn-Wilson conjecture on fp-spectra, which states that the truncated Brown-Peterson spectra generate the category of fp-spectra as a thick subcategory. As a corollary, we deduce the original conjecture at height $1$. Along the way, we prove the existence of $K(n)$-local finite complexes with particularly regular rings of homotopy groups. 
\end{abstract}

\maketitle 

\tableofcontents

\section{Introduction} 

In this paper, we investigate $K(n)$-local spectra whose homotopy groups with coefficients in a type $n$ finite complex are degreewise finite, and characterize them as the thick subcategory generated by the Morava $E$-theory spectrum. As an application of our result, we deduce the height one case of the Hahn-Wilson conjecture. As part of our arguments, we also prove the existence of $K(n)$-local finite complexes of type $n$ whose ring of homotopy groups is an exterior algebra. 

\subsection{Motivation: The Hahn-Wilson conjecture} 

One of the great achievements of chromatic homotopy theory is the classification of finite spectra up to thick subcategories, obtained by Hopkins and Smith as a consequence of the nilpotence theorem \cite{devinatz1988nilpotence, hopkins1998nilpotence}. More precisely, there is a thick subcategory of $p$-local finite spectra for each nonnegative integer $n$, consisting of finite spectra of \emph{type at least $n$}, and this is a complete list of thick subcategories. In practice, finite spectra are extremely difficult from a computational point of view, and to extract useful invariants thereof one takes their homology with respect to certain ``nice''  infinite spectra which can often be understood completely.

In \cite{mahowald1999brown}, Mahowald and Rezk introduce a natural class of such infinite spectra which are particularly amenable to understanding, which they call fp-spectra. We recall that a bounded below, $p$-complete spectrum $X$ is \emph{fp} if it satisfies either of the following equivalent conditions: 
\begin{enumerate}
    \item $\Hrm^{*}(X; \fieldp)$ is finitely presented as a module over the Steenrod algebra, 
    \item there exists a non-zero $p$-local finite spectrum $V$ such that $\pi_{*}(V \otimes X)$ is finite. 
\end{enumerate}
By looking at the thick subcategory of finite $V$ which satisfy the second condition, fp-spectra can be further assigned an \emph{fp-type}, which in this case is an integer $-1 \leq n 
< \infty$. More precisely, if the thick subcategory of finite $V$ satisfying (2) is the category of $p$-local spectra of type at least $n+1$, then we say that $X$ is of fp-type $n$. Many infinite spectra of geometric origin happen to be fp, such as the examples listed in \cref{table:examples_of_fp_types}. 

Mahowald and Rezk show that fp-spectra have many interesting properties; for example, there is a good duality theory in this context, and their chromatic localizations can be calculated using a Tate-like spectral sequence. 

\begin{table}[htbp]
\centering
\caption{Some examples of fp-spectra, implicitly $p$-complete.} 
\begin{tblr}{
  vline{2-7} = {1-2}{},
  hline{2} = {-}{},
}
fp-type  & -1 & 0   & 1   & 2 & \ldots & n  \\
Examples & $\mathbf{Z}/p^{k}$ & $\mathbf{Z}$ & $\mathrm{ko}, \mathrm{ku}$ & $\mathrm{tmf}$ & \ldots  &$\BPn{n}$ 
\end{tblr}
\label{table:examples_of_fp_types}
\end{table}

Another way in which fp-spectra arise naturally is in higher analogues of Quillen-Lichtenbaum conjectures, as conjectured by Ausoni and Rognes \cite{rognes2000algebraic, ausoni2008chromatic}. In the case of the truncated Brown-Peterson spectrum $\BPn{n}$, these were proven in the celebrated work of Hahn and Wilson on chromatic redshift
\cite{hahn2022redshift}, who show that $\K(\BPn{n})^{\wedge}_{p}$ is of fp-type $n+1$.  This highly non-trivial example of an fp-spectrum led them to ask for a general structure result, leading to the following conjecture: 

\begin{conjecture}[Hahn-Wilson\footnote{The authors have learned of this conjecture from Dylan Wilson in 2021, but as far as we are aware, this is the first time it appears in writing. The second author has spoken about this problem during the Homotopy Theory Oberwolfach in 2023.}]
\label{conjecture:introduction_hahn_wilson_conjecture}
The $\infty$-category of spectra of fp-type $\leq n$ is generated as a thick subcategory of $p$-complete spectra by the truncated Brown-Peterson spectrum $\BPn{n}$.\footnote{The statement a priori depends on the form of $\BP\langle n\rangle$, i.e. the choice of generators of $\BP$, but by \cite{leeuniqueness} the underlying spectrum does not depend on these choices at least when $p>2$ or $n\leq 3$. We do not know if this holds when $p = 2$ and $n \geq 4$, and in this case the Hahn-Wilson conjecture should be interpreted as saying that at the very least any form of $\BP \langle n \rangle$ generates the same thick subcategory.}
\end{conjecture}
More informally, the conjecture states that up to thick subcategories, the only fp-spectra are the obvious ones. Before describing our results in this direction, let us mention a few reasons why this statement might be considered plausible: 
\begin{enumerate}
    \item the conjecture is true and easy to verify when $n = -1, 0$, as observed in the original paper of Mahowald and Rezk, 
    \item as part of their analysis of $\K(\BPn{n})^{\wedge}_{p}$ using a proto-version of the even filtration \cite{hahn2022motivic, pstrkagowski2023perfect}, Hahn and Wilson show that it is $\MU$-nilpotent,
    \item if $X$ is an fp-spectrum, then by \cite[{Proposition 3.2}]{mahowald1999brown} the Steenrod algebra module $\Hrm^{*}(X; \fieldp)$ is extended from a module over a finite subalgebra $\acat^\ast(m) \subseteq \acat^{\ast}$ for some finite $m$. 
\end{enumerate}
We are first to admit that arguably all three pieces of evidence are relatively weak, and perhaps the biggest reason we became interested in this problem was our hope that there are yet undiscovered pockets of regularity in the structure of stable homotopy theory. 

Part of the allure of \cref{conjecture:introduction_hahn_wilson_conjecture} is that since $\BP \langle n \rangle$ satisfies the telescope conjecture as an $\MU$-module, it implies the following conjecture from the original paper of Mahowald and Rezk: 

\begin{conjecture}[{Mahowald-Rezk, \cite[{Conjecture 7.3}]{mahowald1999brown}}]
\label{conjecture:introduction_mahowald_rezk}
Let $X$ be of fp-type $n$. Then $X$ satisfies the telescope conjecture in the sense that $L_{n}^{f}X \rightarrow L_{n}X$ is an equivalence. 
\end{conjecture}

Due to the tour de force by Burklund-Hahn-Levy-Schlank \cite{burklund2023k}, it is now known that the telescope conjecture fails for all $n \geq 2$, which arguably was also the prevailing belief at the time of writing of \cite{mahowald1999brown}. In this sense, \cref{conjecture:introduction_mahowald_rezk} is asking about whether the telescope conjecture can fail for ``reasonable'' spectra. Note that the counterexamples given in \cite{burklund2023k} have the property that $\pi_{*}(V \otimes X)$ is bounded for suitable finite $V$, although not finite, so if the Mahowald-Rezk conjecture is true, the line it treads is indeed very thin. 

In the current work, we prove a natural variant of \cref{conjecture:introduction_hahn_wilson_conjecture} obtained, informally, by moving away from the possible failure or truth of the telescope conjecture for fp-spectra. More precisely, we define a natural $K(n)$-local analogue of fp-spectra, and show that they are generated as a thick subcategory by the localization of $\BPn{n}$. Using that the telescope conjecture is true at height one, we are able to use this monochromatic result to deduce the following case of \cref{conjecture:introduction_hahn_wilson_conjecture}: 

\begin{theorem}[{\cref{theorem:main_text_hahn_wilson_conjecture_at_height_one}}]
\label{theorem:introduction_hahn_wilson_at_height_one}
The Hahn-Wilson conjecture is true when $n = 1$; that is, spectra of fp-type $\leq 1$ coincide with the thick subcategory generated by $\BPn{1}$. 
\end{theorem}

We remark that the deduction of 
\cref{theorem:introduction_hahn_wilson_at_height_one} from the monochromatic statement is not particularly difficult. The majority of the current work is devoted to proving the $K(n)$-local result, which we describe now. 

\subsection{Monochromatic fp-spectra} 

Since the Eilenberg-MacLane spectrum $\fieldp$ vanishes $K(n)$-locally, to generalize the notion of being fp to the $K(n)$-local setting, it is natural to use conditions on homotopy groups instead. 

\begin{definition}
\label{definition:introduction_local_fp_type}
We say a $K(n)$-local spectrum $X$ is \emph{locally fp} if $\pi_{*}(V \otimes X)$ is degreewise finite for any finite complex $V$ of type $n$. 
\end{definition}

This definition is motivated by the following easy observation:

\begin{example}
Let $X$ be a spectrum of fp-type $n$ and assume that it satisfies the telescope conjecture at height $n$ in the sense that $L_{T(n)} X \simeq L_{K(n)} X$. Then $L_{K(n)} X$ is locally fp, as we show in \cref{lemma:spectra_of_fp_type_n_which_satisfy_the_telescope_localize_to_local_fps}. 
\end{example}

As observed in \cite{mahowald1999brown}, chromatic localizations of finite spectra can often be written as a chromatic localization of an fp-spectrum. In line with this observation, we have the following:  

\begin{example}
Let $X$ be a $K(n)$-locally dualizable spectrum; for example, a $K(n)$-localization of a finite spectrum. Then $X$ is locally fp by \cite[{Theorem 8.5, Corollary 8.12}]{hovey1999morava}. 
\end{example}

It is not difficult to see that not all locally fp-spectra are $K(n)$-locally dualizable. For example, $L_{K(n)} \BPn{n}$ is locally fp but not dualizable. The following theorem, which is the main result of the current work, shows that this is essentially the only new example: 

\begin{theorem}
\label{theorem:introduction_easy_statement_of_kn_local_hahn_wilson}
Locally fp-spectra are exactly the thick subcategory of $K(n)$-local spectra generated by $L_{K(n)} \BPn{n}$. 
\end{theorem}

We stated our result in the above way to make the analogy with \cref{conjecture:introduction_hahn_wilson_conjecture} clearer. However, since we are working $K(n)$-locally, it is arguably more natural to state the result in terms of Morava $E$-theory, as we do now. 

We write $E$ for the Morava $E$-theory spectrum associated to the Honda formal group law over $\mathbf{F}_{p^{n}}$. As a spectrum, it is equivalent to a finite direct sum of $L_{K(n)} \BPn{n}$\footnote{The $K(n)$-localization of $\BPn{n}$ can be identified with what is classically known as the completed Johnson-Wilson theory $\widehat{E(n)}$ by \cite[{Theorem 4.1}]{baker1989liftings}. There is a map of rings $\widehat{E(n)}_{*} \rightarrow E(\mathbf F_{p^n}, \mathbf K)_{*}$ (compatible with formal groups on both sides) which presents the target as a finite free module over the source. Here, $E(\mathbf F_{p^n}, \mathbf K)$ is the Lubin-Tate spectrum associated with some formal group $\mathbf K$ over $\mathbf F_{p^n}$. This lifts to a map of homotopy ring spectra by the Landweber exact theorem, and together with \cite{luecke2022arentmoravaetheories}, this gives us the needed result.}; in particular, they generate the same thick subcategory. We deduce \cref{theorem:introduction_easy_statement_of_kn_local_hahn_wilson} as a consequence of the following more elaborate characterization: 

\begin{theorem}
[{\cref{theorem:classification_of_local_fp_types_in_terms_of_cohomology_and_morava_e_theory}}]
\label{theorem:introduction_local_fp_types_in_terms_of_morava_e_theory}
For a $K(n)$-local spectrum $X$, the following are equivalent: 

\begin{enumerate}
    \item $X$ is locally fp in the sense of \cref{definition:introduction_local_fp_type}, 
    \item $E^{*}X$ is finitely generated as a module over the ring of operations $E^{*}E$,
    \item $X$ belongs to the thick subcategory generated by $E$, 
    \item $X$ is reflexive; that is, the canonical map $X \rightarrow D^{2}(X)$ into its $K(n)$-local monoidal double dual is an equivalence.
\end{enumerate}
\end{theorem}
At height one, which is already sufficient to deduce \cref{theorem:introduction_hahn_wilson_at_height_one}, the above result was proven by Rebekah Hahn and Steve Mitchell in their study of Iwasawa theory of $K(1)$-local spectra \cite[{Theorem 8.1}]{hahn2007iwasawa}. 

We observe that the second condition in \cref{theorem:introduction_local_fp_types_in_terms_of_morava_e_theory} is analogous to the characterization of integral fp-spectra in terms of their mod $p$-cohomology. The fourth characterization, which was suggested to us by Strickland, is interesting because it describes locally fp-spectra purely in terms of the monoidal structure. This variety of different descriptions gives further evidence that locally fp-spectra form a very natural class of spectra to consider. 

We now briefly sketch the arguments that go into the proof of \cref{theorem:introduction_local_fp_types_in_terms_of_morava_e_theory}. The implication $(2) \Rightarrow (3)$ is relatively standard. The ring of operations can be identified with the twisted group algebra $E^{*}E \simeq E^{*} \llbracket \mathbf{G}_{n} \rrbracket$ associated to the Morava stabilizer group $\mathbf{G}_{n}$ acting on the $2$-periodic Lubin-Tate ring. The group $\mathbf{G}_{n}$ is $p$-adic analytic, and so as a consequence of the work of Lazard it is virtually of finite cohomological dimension. 

We use the work of Lazard to show that $E^{*}E$ is a finite algebra over a ring of finite global dimension, see \cref{proposition:twisted_group_algebra_of_a_uniform_subgroup_is_finite_dim_and_notherian}. The implication $(2) \Rightarrow (3)$ is proven by lifting a finite projective resolution of $E^{*}X$ to a resolution of $X$ itself, and a reduction argument allowing one to work over a finite Galois extension of the $K(n)$-local sphere, see \cref{lemma:being_in_thick_subcategory_of_e_up_to_finite_galois_extension_means_youre_in_thick_subcat_of_e}. 

$(3) \Rightarrow (4)$ follows from a result of Strickland that $E$ itself is reflexive. To prove $(4) \Rightarrow (1)$, we use the related fact that the $K(n)$-local Brown-Comenetz dual spectrum is invertible. 

The hard part of \cref{theorem:introduction_local_fp_types_in_terms_of_morava_e_theory} is proving that $(1) \Rightarrow (2)$. While $E^{*}X$ and $\pi_{*}(X)$ are related via the $K(n)$-local Adams-Novikov spectral sequence, we were unable to make an analysis of this spectral sequence work at heights $n > 1$. Instead, we proceed differently, by studying homotopy groups of certain finite complexes, which we describe next.

\subsection{Homotopy of \texorpdfstring{$K(n)$}{K(n)}-local Moore spectra} 

Our strategy to prove that $E^{*}(X)$ is finitely generated over $E^{*}E$ when $X$ is locally fp is based on the construction of a particularly nice $K(n)$-local finite complex $V$ of type $n$. We begin by considering a generalized Moore complex
\[
M \colonequals \thesphere / (p^{\alpha_{0}}, v_{1}^{\alpha_{1}}, \ldots, v_{n-1}^{\alpha_{n-1}}). 
\]
By the work of Burklund, for large enough indices $\alpha_{i}$, $M$ admits a structure of an $\mathbf{E}_{2}$-ring spectrum \cite{burklund2022multiplicative}. We have $E_{*}(M) \simeq E_{*}/I$ with $I = (p^{\alpha_{0}}, \ldots, v_{n-1}^{\alpha_{n-1}})$.  

For a suitably chosen open subgroup $U \leq \mathbf{G}_{n}$, we write $V \colonequals E^{h U} \otimes M$; this is the desired $K(n)$-local finite complex of type $n$. The inclusion of homotopy fixed points provides a map of $\mathbf{E}_{2}$-rings 
\[
V \rightarrow E \otimes M. 
\]
Our aim is to analyze the universal coefficient spectral sequence
\[
\Ext_{V_{*}}(V_{*}(X), E_{*}M) \Rightarrow (E \otimes M)^{*}X
\]
associated to this map, deducing information about $(E \otimes M)^{*}(X)$ from that of $V_{*}(X)$, which by assumption that $X$ is locally fp is degreewise finite. Finally, information about $E^{*}X$ itself is obtained by a completeness argument. 

To make this work, one wants $V_{*}$ to be a well-behaved ring, which arguably is a tall order for homotopy groups of a finite complex, even $K(n)$-locally. The following surprising result shows, however, that such well-behaved finite complexes do exist: 

\begin{theorem}[{\cref{theorem:existence_of_a_convenient_generalized_moore_and_open_subgroup}}]
\label{theorem:introduction_existence_of_nice_finite_complexes}
For appropriately chosen indices $(\alpha_{0}, \ldots, \alpha_{n-1})$ and an appropriately chosen normal open subgroup $U \leq \mathbf{G}_{n}$, the $K(n)$-local finite complex $V = E^{hU} \otimes M$ has the following properties: 
\begin{enumerate}
\item $U$ acts trivially on $E_{*}(M) \simeq E_{*}/I$ and 
\[
\Hrm^{*}(U, E_{*}/I) \simeq \Lambda_{E_{*}/I}(x_{1}, \ldots, x_{d}) 
\]
is an exterior algebra over $E_{*}/I$ on any basis $x_{i} \in \Hrm^{1}(U, \ZZ/p^{\alpha_{0}})$, where $d = n^{2}$ is the dimension of the Morava stabilizer group, 
\item the $K(n)$-local Adams-Novikov spectral sequence 
\[
\Hrm^{*}(U, E_{*}/I) \Rightarrow \pi_{\ast}(V)
\]
associated with $V = E^{hU} \otimes M$ collapses on the second page, 
\item any collection of lifts $\overline{x}_{i} \in \pi_{-1}(V)$ of $x_{i}$ induces an isomorphism of rings
\[
\pi_{\ast}(V) \simeq \Lambda_{E_{*}/I}(\overline{x}_{1}, \ldots, \overline{x}_{d}).
\]
\end{enumerate}
\end{theorem}

We remark here that the indices $\alpha_{i}$ can be chosen to be arbitrarily large; in particular, using the work of Burklund one can make $V$ into an $\mathbf{E}_{k}$-algebra for any finite $k$. Moreover, while the choice of appropriate $U$ depends on these indices, the subgroup can be chosen to be arbitrarily small. Note that since the subgroup is open, $E^{hU}$ is a finite Galois extension of the $K(n)$-local sphere, so one way to interpret \cref{theorem:introduction_existence_of_nice_finite_complexes} is that up to finite Galois extensions, homotopy rings of generalized Moore complexes have a particularly regular structure. 

The proof of \cref{theorem:introduction_existence_of_nice_finite_complexes} is surprisingly delicate, and it involves some structural properties of the lower $p$-series of a $p$-adic analytic group, as well as a simultaneous analysis of the $K(n)$-local Adams-Novikov spectral sequences $\Hrm^{*}(U, E_{*}/I) \Rightarrow \pi_{\ast}(V)$ where both the indices $\alpha_{i}$ and the subgroup $U$ are allowed to vary. 

We finish the introduction by sketching the proof of part $(1) \Rightarrow (2)$ of \cref{theorem:introduction_local_fp_types_in_terms_of_morava_e_theory} using the $K(n)$-local finite complex $V$ of \cref{theorem:introduction_existence_of_nice_finite_complexes}. The key observation is that the second page of the universal coefficient spectral sequence
\begin{equation}
\label{equation:introduction:collapsing_uct_cohomology_sseq}
\Ext_{V_{*}}(E_{*}/I, E_{*}/I) \Rightarrow \pi_{*}(\map_{\Mod_{V}}(E \otimes M, E \otimes M)) 
\end{equation}
is given by the cohomology of an exterior algebra over the graded Artinian ring $E_{*}/I$, and so is given by a polynomial algebra on classes of even total degree. In particular, it is noetherian and the spectral sequence collapses on the second page.

Now we assume we have a locally fp-spectrum $X$, so that $V_{*}(X)$ is degreewise finite. We show that this implies that it is finitely presented over $V_{*}$. Using noetherianness, we use this to prove that the second page of the universal coefficient spectral sequence 
\begin{equation}
\label{equation:introduction_sseq_approximating_cohomology_of_a_local_fp_type}
\Ext_{V_{*}}(V_{*}(X), E_{*}/I) \Rightarrow \pi_{*}(\map_{\Mod_{V}}(V \otimes X, E \otimes M)) \simeq (E \otimes M)^{-*}(X) 
\end{equation}
is a finitely presented module over the second page of (\ref{equation:introduction:collapsing_uct_cohomology_sseq}). Then a classical argument involving noetherianness shows that (\ref{equation:introduction_sseq_approximating_cohomology_of_a_local_fp_type}) collapses at a finite page, which is then also finitely presented. Using this, we deduce that $(E \otimes M)^{*}X$, and then $E^{*}X$ itself, is finitely presented, ending the argument. 

\subsection{Acknowledgements} 

We would like to thank Dustin Clausen, Paul Goerss, Jeremy Hahn, Ishan Levy, John Rognes, Andy Senger, Neil Strickland, and Dylan Wilson for helpful discussions related to this project. The second author would like to thank the many Starbucks of Osaka\footnote{In particular, the ones at \begin{CJK}{UTF8}{min}大阪市中央区今橋２丁目２−２\end{CJK} and \begin{CJK}{UTF8}{min}大阪市中央区難波５丁目１−６０\end{CJK}.} where this paper was completed. 

\section{Algebraic preliminaries} 

In this short section, we collect algebraic results regarding cohomology of $p$-adic analytic groups and the twisted group rings associated to subgroups of the Morava stabilizer group. 

\subsection{Cohomology of uniform groups} 

We recall that a profinite group $U$ is said to be \emph{powerful} if it is pro-$p$ and the quotient $U/\overline{U^{p}}$ by the closed subgroup generated by $p$-th powers is abelian (resp. $U/\overline{U^{4}}$ is abelian if $p = 2$). It is said to be uniformly powerful (or \emph{uniform}) if it is in addition finitely generated and torsion-free; see \cite{dixon2003analytic} for a comprehensive treatment. 

A celebrated result of Lazard shows that if $U$ is uniform, then its cohomology algebra
\[
\Hrm^{*}(U, \fieldp) \simeq \Lambda_{\fieldp}(x_{1}, \ldots, x_{d}) 
\]
is exterior on classes $x_{i} \in \Hrm^{1}(U, \fieldp)$. In this short section, we consider a natural strengthening of the notion of being uniform, and show that under this assumption the isomorphism of Lazard can be lifted to one with coefficients in $\ZZ/p^{n}$. 

\begin{notation}
\label{notation:lower_p_series}
If $U$ is a pro-$p$ group, we write 
\[
U = U_{1} \geq U_{2} \geq U_{3} \geq \ldots 
\]
for its lower $p$-series defined inductively by $U_{1} \colonequals U$ and $U_{k+1} \colonequals \overline{U_{k}^{p} \cdot [U_{k}, U]}$. This is a sequence of topologically characteristic closed subgroups. 
\end{notation}

\begin{remark}
We will only be interested in the above notion when $U$ is uniform, in which case \cite[{Theorem 3.6}]{dixon2003analytic} shows that we have an equality of subsets
\[
U_{k+1} = \{ u^{p^{k}} \ | \ u \in U \}. 
\]
In particular, in this case the subset of $p^{k}$-th powers is already closed and a subgroup. 
\end{remark}

\begin{definition}
\label{definition:n_abelian_uniforn_group}
Let $U$ be a uniform pro-$p$ group and let $n \geq 1$. We say that $U$ is \emph{$n$-uniform} if 
\begin{enumerate}
\item $U / U_{n+1}$ is abelian and $p > 2$, 
\item $U / U_{n+2}$ is abelian and $p = 2$. 
\end{enumerate}
\end{definition}

\begin{example}
A uniform group is $1$-uniform, by definition. 
\end{example}

\begin{example}
A uniform group is abelian if and only if it is $n$-uniform for all finite $n$. By a structure result for uniform groups \cite[{Theorem 4.17}]{dixon2003analytic}, the only examples are $\ZZ_{p}^{\oplus d}$. 
\end{example}

\begin{lemma}
\label{lemma:subgroup_of_pth_powers_of_an_n_abelian_gp_is_nplusone_abelian}
Let $U$ be an $n$-uniform group. Then $U_{2} = \{ u^{p} \ | \ u \in U \}$ is $(n+1)$-uniform. 
\end{lemma}

\begin{proof}
We will make use of the commutator property of the lower $p$-series, which states that 
\[
[U_{k}, U_{l}] \leq U_{k+l},
\]
see \cite[Proposition 1.16, (b)]{dixon2003analytic}. We only do the case of $p > 2$, since $p = 2$ differs only by a shift by one in indexing. 

As we have $(U_{2})_{n+2} = U_{n+3}$, the claim is equivalent to showing that $[U_{2}, U_{2}] \leq U_{n+3}$. Since $U$ is $n$-uniform, we know that $[U_{1}, U_{1}] \leq U_{n+1}$. By the commutator property, $U_{n+1}/U_{n+2}$ is central in $U/U_{n+2}$. Since the commutator is linear in each variable up to conjugation, see \cite[Prelude, 0.2.(i)]{dixon2003analytic}, centrality implies that 
\[
[u^{p}, v] \equiv [u, v]^{p} \mod U_{n+2}. 
\]
Since any element of $U_{2}$ can be written as a $p$-th power, this implies that $[U_{2}, U_{1}] \leq U_{n+2}$. Since $U_{n+2}/U_{n+3} \leq U/U_{n+3}$ is also central, the same argument shows that $[U_{2}, U_{2}] \leq U_{n+3}$. 
\end{proof}

\begin{proposition}
\label{proposition:zpn_cohomology_of_n_abelian_uniform_group_is_exterior}
Let $U$ be an $n$-uniform group. Then 
\begin{enumerate}
\item $\Hrm^{1}(U, \ZZ/p^{n})$ is free over $\ZZ/p^{n}$ of rank equal to the dimension $d = \mathrm{dim}(U)$, 
\item any choice of a basis $x_{1}, \ldots, x_{d} \in \Hrm^{1}(U, \ZZ/p^{n})$ induces an isomorphism 
\[
\Hrm^{\ast}(U, \ZZ/p^{n}) \simeq \Lambda_{\ZZ/p^{n}}(x_{1}, \ldots, x_{d}) 
\]
between the cohomology and an exterior algebra. 
\end{enumerate}
\end{proposition}

\begin{proof}
If $n = 1$, this is a classical result due to Lazard \cite[{Theorem 5.1.5}]{symonds2000cohomology}. If $n > 1$, then 
\[
U/ U_{n+1}
\]
is an abelian group of order $p^{nd}$ and rank $d$, so that 
\[
U/U_{n+1} \simeq (\ZZ/p^{n})^{\times d}.
\]
It follows that the cohomology group can be described explicitly as 
\[
\Hrm^{1}(U, \ZZ/p^{n}) \simeq \mathrm{Hom}(U, \ZZ/p^{n}) \simeq \mathrm{Hom}(U/U_{n+1}, \ZZ/p^{n}) \simeq  (\ZZ/p^{n})^{\times d},
\]
so that the reduction map 
\[
\mathrm{red} \colon \Hrm^{1}(U, \ZZ/p^{n}) \rightarrow \Hrm^{1}(U, \ZZ/p)
\]
is surjective. Since $\ZZ/p$-cohomology algebra is generated in degree one, we deduce that 
\[
\mathrm{red} \colon \Hrm^{\ast}(U, \ZZ/p^{n}) \rightarrow \Hrm^{\ast}(U, \ZZ/p)
\]
is surjective in all degrees. It follows that each $\ZZ/p^{n}$-cohomology group is free over $\ZZ/p^{n}$ of rank the same as the $\ZZ/p$-dimension of the corresponding $\ZZ/p$-cohomology group. 

Let $x_{1}, \ldots, x_{d} \in \Hrm^{1}(U, \ZZ/p^{n})$ be a basis as in the statement. We claim that $x_{i}^{2} = 0$. If $p > 2$, this follows immediately from commutativity. If $p = 2$, then the stricter definition of being $n$-uniform at the even prime together with the above argument shows that $\ZZ/2^{n+1}$-cohomology is also degreewise free and that 
\[
\Hrm^{\ast}(U, \ZZ/2^{n+1}) \rightarrow \Hrm^{\ast}(U, \ZZ/2^{n})
\]
is also surjective. Thus, $x_{i}$ lift to some $\ZZ/2^{n+1}$-cohomology elements $\overline{x}_{i}$. By anticommutativity, $\overline{x}_{i}^{2}$ are $2$-torsion, so that their images $x_{i}^2$ in $\ZZ/2^{n}$-cohomology are zero. 

Finally, let $x_{1}, \ldots, x_{d} \in \Hrm^{1}(U, \ZZ/p^{n})$ be a basis as in the statement. Since $x_{i}^{2} = 0$, we have an induced map of algebras
\[
\Lambda_{\ZZ/p^{n}}(x_{1}, \ldots, x_{d}) \rightarrow \Hrm^{*}(U, \ZZ/p^{n}).
\]
This is a map of free $\ZZ/p^{n}$-modules of the same rank which is an isomorphism mod $p$ and hence is an isomorphism, as needed. 
\end{proof}

\subsection{Twisted group algebras} 

Associated to a profinite group $G$ we have the completed group algebra 
\[
\fieldp \llbracket G \rrbracket \colonequals \varprojlim \fieldp [G/N],
\]
where the limit is taken over open normal subgroups $N \lhd G$. A classical result shows that if $G$ is uniform, then the associated graded of $\fieldp \llbracket G \rrbracket$ with respect to the powers of the augmentation ideal is isomorphic to a polynomial algebra \cite[{Theorem 7.22}]{dixon2003analytic}. Using methods of filtered ring theory \cite{huishi}, one deduces that the complete group algebra itself is left (and right) noetherian and of global dimension equal to the dimension of $G$. 

In this short section, we show that analogous properties hold for the twisted group algebras associated to the actions of the Morava stabilizer $\mathbf{G}_{n}$ group on coefficients of Morava $E$-theory. The key property is that $\mathbf{G}_{n}$ is $p$-adic analytic; that is, it has a uniformly powerful open subgroup. For more information about this action, together with the notation we use, see the excellent overview of Barthel-Beaudry \cite[{\S 3.1}]{barthel2020chromatic}. 

\begin{recollection}
We recall the construction of the needed group algebra from \cite[{\S 5}]{hovey2004operations}. For each $k \geq 0$, the quotient $E^{*}/\mfrak^{k}$ has the discrete topology. If $U$ is an open subgroup of the Morava stabilizer group, we define 
\[
(E^{*} / \mfrak^{k}) \llbracket U \rrbracket \colonequals \varprojlim (E^{*} / \mfrak^{k})[U/V],
\]
the limit taken over the poset of those open normal $V \leq U$ such that the action on $E^{*}/\mfrak^{k}$ factors through $U/V$. The completed twisted group algebra of $U$ is given by the sequential limit 
\[
E^{*} \llbracket U \rrbracket \colonequals \varprojlim (E^{*} / \mfrak^{k}) \llbracket U \rrbracket.
\]
\end{recollection}

The importance of completed group algebras from our perspective comes from the natural way they arise as endomorphisms of Morava $E$-theory. The following is presumably well-known to experts but we could not find a suitable reference at this level of generality: 

\begin{lemma}
\label{lemma:ehu_linear_endo_of_e_is_a_twisted_group_algebra}
Let $U \leq \mathbf{G}_{n}$ be an open subgroup. Then the action of $U$ on $E$ induces an isomorphism of graded rings 
\begin{equation}
\label{equation:ehu_linear_endomorphisms_of_e_as_a_group_algebra}
E^{*} \llbracket U \rrbracket \simeq \pi_{-*} \map_{\Mod_{E^{hU}}}(E, E).
\end{equation}
In particular, $E^{*}E \simeq E^{*} \llbracket \mathbf{G}_{n} \rrbracket$. 
\end{lemma}

\begin{proof}
If $U = \mathbf{G}_{n}$, this is a result of Hovey \cite[{Theorem 5.5}]{hovey2004operations}.  To prove the general case, it will be convenient to think of the continuous action of $\mathbf{G}_{n}$ on Morava $E$-theory as defining a sheaf
\[
S \mapsto E(S) 
\]
of $K(n)$-local spectra on profinite $\mathbf{G}_{n}$-sets as in \cite[{\S 3.1}]{mor2023picard}. In particular, if $K \leq \mathbf{G}_{n}$ is a closed subgroup, we have $E(K) = E^{hK}$, the fixed points of Devinatz and Hopkins \cite{devinatz2004homotopy}.

We can equivalently show the continuous linear dual of (\ref{equation:ehu_linear_endomorphisms_of_e_as_a_group_algebra}), namely that there's an isomorphism of $E_{*}$-modules
\[
\pi_{*}(E \widehat{\otimes}_{E^{hU}} E) \simeq \map_{\mathrm{cts}}(U, E_{*}),
\]
where $\widehat{\otimes}$ denotes the $K(n)$-local tensor product, which we consider as an $E$-module using the left factor. The relative tensor product on the left hand side can be calculated as the $K(n)$-local colimit of the simplicial diagram
\[
\begin{tikzcd}[column sep=1.5em]
\cdots
  \arrow[r, shift left=3pt]
  \arrow[r]
  \arrow[r, shift right=3pt]
&
E \widehat{\otimes} E^{hU} \widehat{\otimes} E
  \arrow[r, shift left=2pt]
  \arrow[r, shift right=2pt]
&
E \widehat{\otimes} E
\end{tikzcd}
\]
which using the notation of the first paragraph we can identify with the diagram 
\begin{equation}
\label{equation:diagram_defining_relative_tensor_e_ehu_e_using_profinite_gn_sets}
\begin{tikzcd}[column sep=1.5em]
\cdots
  \arrow[r, shift left=3pt]
  \arrow[r]
  \arrow[r, shift right=3pt]
&
E(\mathbf{G}_{n} \times \mathbf{G}_{n}/U \times \mathbf{G}_{n})
  \arrow[r, shift left=2pt]
  \arrow[r, shift right=2pt]
&
E(\mathbf{G}_{n} \times \mathbf{G}_{n})
\end{tikzcd}
\end{equation}
induced by the underlying cosimplicial object of the augmented cosimplicial $\mathbf{G}_{n}$-set 
\[
\begin{tikzcd}[column sep=1.5em]
\mathbf{G}_{n} \times_{\mathbf{G}_{n}/U} \mathbf{G}_{n}
  \arrow[r]
&
\mathbf{G}_{n} \times \mathbf{G}_{n}
  \arrow[r, shift left=2pt]
  \arrow[r, shift right=2pt]
&
\mathbf{G}_{n} \times \mathbf{G}_{n}/U \times \mathbf{G}_{n}
  \arrow[r, shift left=3pt]
  \arrow[r]
  \arrow[r, shift right=3pt]
&
\cdots
\end{tikzcd}
\]
The latter induces an augmentation of (\ref{equation:diagram_defining_relative_tensor_e_ehu_e_using_profinite_gn_sets}) 
which makes it into a colimit diagram, as can be verified by taking completed $E$-homology, so that 
\[
E \widehat{\otimes} _{E^{hU}} E \simeq E(\mathbf{G}_{n} \times_{\mathbf{G}_{n}/U} \mathbf{G}_{n}).
\]
We deduce that 
\[
\pi_{*}(E \widehat{\otimes} _{E^{hU}} E) \simeq \pi_{*}(E(\mathbf{G}_{n} \times_{\mathbf{G}_{n}/U} \mathbf{G}_{n})) \simeq \map_{\mathrm{cts}}(\mathbf{G}_{n} \times_{\mathbf{G}_{n}/U} \mathbf{G}_{n}, E_{*})^{h\mathbf{G}_{n}},
\]
where the second isomorphism uses that $\mathbf{G}_{n} \times_{\mathbf{G}_{n}/U} \mathbf{G}_{n}$ is a free $\mathbf{G}_{n}$-set. Furthermore, we have
\[
\map_{\mathrm{cts}}(\mathbf{G}_{n} \times_{\mathbf{G}_{n}/U} \mathbf{G}_{n}, E_{*})^{\mathbf{G}_{n}} \simeq \map_{\mathrm{cts}}(U, E_{*}), 
\]
where the isomorphism is given by restricting a $\mathbf{G}_{n}$-equivariant function $\mathbf{G}_{n} \times_{\mathbf{G}_{n}/U} \mathbf{G}_{n} \rightarrow E_{*}$ along the inclusion $U \simeq \ast \times_{\ast} U \hookrightarrow \mathbf{G}_{n} \times_{\mathbf{G}_{n}/U} \mathbf{G}_{n}$. This ends the proof. 
\end{proof}

\begin{lemma}
\label{lemma:global_dimension_of_rings_complete_wrt_to_a_central_element}
Let $R$ be a ring together with an element $x \in R$ such that: 
\begin{enumerate}
    \item $x$ is central, 
    \item $x$ is a non-zero divisor, 
    \item $R$ is $(x)$-complete,
    \item $R/x$ is left noetherian, right noetherian and of finite global dimension $d$. 
\end{enumerate}
Then $R$ is left noetherian, right noetherian and of global dimension $d+1$. 
\end{lemma}

\begin{proof}
This is a combination of \cite[{Chapter I, Proposition 7.1.2, Corollary 7.2.1, Chapter II, Theorem 5.1.5}]{huishi}. 
\end{proof}

The following lemma is useful in reducing from a twisted to an untwisted group algebra: 

\begin{lemma}
\label{lemma:untwisting_twisted_group_algebras}
Let $k \rightarrow k'$ be a finite Galois extension with Galois group $G$, $P$ a profinite group acting on $k'$ through a surjective group homomorphism $P \rightarrow G$ and let $k' \llbracket P \rrbracket$ be the associated completed twisted group algebra. Then, a choice of a section $G \rightarrow P$ (as a map of sets) induces an isomorphism
\[
k' \otimes_{k} k'\llbracket P \rrbracket \simeq \mathrm{M}_{|G|}(k' \llbracket \mathrm{ker}(P \rightarrow G) \rrbracket)
\]
of $k'$-algebras between the base-change of $k' \llbracket P \rrbracket$ and the matrix algebra of the group algebra associated to the kernel. 
\end{lemma}

\begin{proof}
Let $U$ be an open normal subgroup of $P$ contained in $\mathrm{ker}(P\to G)$ and $P' = P/U$ the associated finite quotient. By Galois theory, $k' \otimes_{k} k' \simeq \map(G, k')$, where on the right $G$ acts by permutation, so that 
\[
k' \otimes_{k} k'[P'] \simeq \map(G, k')[P']. 
\]
The right hand side can be identified with the $k'$-linear groupoid algebra associated to the action groupoid of $P'$ acting on $G$ by permutation. The latter is a connected groupoid with $|G|$ objects and automorphism group isomorphic to $\mathrm{ker}(P' \rightarrow G)$. Thus, \cite[{Theorem 3.2}]{Steinberg2006Mobius} associates to a choice of a section $G \rightarrow P'$ an isomorphism 
\[
k' \otimes_{k} k'[P'] \simeq \mathrm{M}_{|G|}(k'[\mathrm{ker}(P' \rightarrow G)]).
\]
Passing to the limit over finite quotients yields the required isomorphism. 
\end{proof}

\begin{proposition}
\label{proposition:twisted_group_algebra_of_a_uniform_subgroup_is_finite_dim_and_notherian}
Let $U \leq \mathbf{G}_{n}$ be an open subgroup of the Morava stabilizer group which is torsion-free. Then the graded ring $E^{*} \llbracket U \rrbracket$ is left (resp. right) noetherian and of left (resp. right) global dimension $n^{2}+n$ in the graded sense. 
\end{proposition}

\begin{proof}
As a twisted group algebra, the ring is isomorphic to its opposite, so we can focus on the category of left modules. Since $E^{*} \llbracket U \rrbracket$ is concentrated in even degrees, its category of left modules is equivalent to the product of the categories of left modules concentrated in even degrees and in odd degrees, which are isomorphic to each other. Thus, it is enough to show that the category $\Mod_{E^{*} \llbracket U \rrbracket}(\mathrm{grAb}^{ev})$ of left modules concentrated in even degrees is noetherian and of global dimension $n^{2}+n$. 

The sequence $p, v_{1}, \ldots, v_{n-1}$ is a sequence of elements to which \cref{lemma:global_dimension_of_rings_complete_wrt_to_a_central_element} can be applied inductively (as $v_{i}$ is fixed by the Morava stabilizer group modulo $(p, v_{1}, \ldots, v_{i-1})$), so the graded analogue of the latter shows that it is enough to show that the category of even graded left modules over 
\[
E^{*}/(p, \ldots, v_{n-1}) \llbracket U \rrbracket \simeq \mathbf{F}_{p^{n}}[u^{\pm 1}] \llbracket U \rrbracket 
\]
is noetherian and of global dimension $n^{2}$. Since this ring has an invertible variable $u$ of degree $2$, it is strongly evenly graded and thus by \cite[{Theorem 2.8}]{dade1980group} the restriction of a module to its degree zero part determines an equivalence of categories with \emph{ungraded} modules over $\mathbf{F}_{p^{n}} \llbracket U \rrbracket$. 

We consider the composite $U \hookrightarrow \mathbf{G}_{n} \rightarrow \mathrm{Gal}(\mathbf{F}_{p^{n}} / \mathbf{F}_{p})$ and write $Q$ for its image and $K$ for its kernel, so that there is a short exact sequence of profinite groups
\[
0 \rightarrow K \rightarrow U \rightarrow Q \rightarrow 0.
\]
We write $q \colonequals |Q|$, so that $\mathbf{F}_{p^{n/q}} = (\mathbf{F}_{p^{n}})^{Q}$ and we have an identification $Q \simeq \mathrm{Gal}(\mathbf{F}_{p^{n}} / \mathbf{F}_{p^{n/q}})$. By construction, $\mathbf{F}_{p^{n/q}}$ is central in the group algebra and by \cref{lemma:untwisting_twisted_group_algebras} there exists an isomorphism
\[
\mathbf{F}_{p^{n}} \otimes_{\mathbf{F}_{p^{n/q}}} \mathbf{F}_{p^{n}} \llbracket U \rrbracket \simeq \mathrm{M}_{q}(\mathbf{F}_{p^{n}} \llbracket K \rrbracket)
\]
of $\mathbf{F}_{p^{n}}$-algebras. Since passing to a finite separable extension does not change the global dimension or the property of being left noetherian
\cite[{Theorem 5.1}]{Kadison1995SplitSeparable},
by Morita invariance we deduce that it is enough to show that
$\mathbf{F}_{p^{n}} \llbracket K \rrbracket$ is left noetherian and of global dimension $n^{2}$ .
Since $K$ acts trivially on $\mathbf{F}_{p^{n}}$, this is the ordinary (untwisted) group algebra.

As $K$ is torsion-free and $p$-adic analytic of dimension $n^{2}$, $\mathbf{F}_{p^{n}} \llbracket K \rrbracket \simeq \mathbf{F}_{p^{n}} \otimes_{\mathbf{F}_{p}} \mathbf{F}_{p} \llbracket K \rrbracket$ has the needed properties by a combination of \cite[{Theorem J}]{ardakov_brown_primeness_2027} and \cite[{Corollary 7.25}]{dixon2003analytic}.
\end{proof}

\begin{corollary}
\label{corollary:twisted_group_algebra_of_gn_is_noether}
For any open subgroup $U \leq \mathbf{G}_{n}$, the ring $E^{*} \llbracket U \rrbracket$ is left and right noetherian. In particular, $E^{*}E$ is left and right noetherian. 
\end{corollary}

\begin{proof}
Since $U$ is open, it contains an open uniform subgroup $U' \leq U$. The needed statement then follows from \cref{proposition:twisted_group_algebra_of_a_uniform_subgroup_is_finite_dim_and_notherian}, since the map of rings 
\[
E^{\ast} \llbracket U' \rrbracket \hookrightarrow E^{\ast} \llbracket U \rrbracket 
\]
presents the target as a free left (or right) module on basis given by any set of representatives for the finite quotient $U / U'$. Since the source is left and right noetherian, we deduce that so is the target. 
\end{proof}

\begin{proposition}
\label{proposition:twisted_group_algebra_of_a_pro_p_group_is_local}
Let $U \leq \mathbf{G}_{n}$ be an open subgroup of the Morava stabilizer group which is a pro-$p$-group and which lies in the kernel of $\mathbf{G}_{n} \rightarrow \mathrm{Gal}(\mathbf{F}_{p^{n}} / \mathbf{F}_{p})$. Then $E^{0} \llbracket U \rrbracket$ is a local ring. 
\end{proposition}

\begin{proof}
We have to show that the Jacobson radical is maximal. Let $\mfrak = (p, \ldots, u_{n-1}) \subseteq E^{0} \llbracket U \rrbracket$ be the double-sided ideal generated by the maximal ideal of $E^{0}$. By construction as a limit, $E^{0} \llbracket U \rrbracket$ is $\mfrak$-complete, so that $\mfrak$ is contained in the Jacobson radical. It follows that it is enough to show that 
\[
(E^{0} \llbracket U \rrbracket) / (p, \ldots, u_{n-1}) \simeq \mathbf{F}_{p^{n}} \llbracket U \rrbracket 
\]
is local. Note that by the assumption on $U$, the right hand side is actually an ordinary (untwisted) group algebra. In particular, augmentation provides a map of rings 
\[
\mathbf{F}_{p^{n}} \llbracket U \rrbracket \rightarrow \mathbf{F}_{p^{n}}. 
\]
As $U$ is pro-$p$, the source is complete with respect to filtration by powers of the augmentation ideal by \cite[{\S 7.1}]{dixon2003analytic} and thus the kernel is contained in the Jacobson radical. This ends the argument. 
\end{proof}

\begin{corollary}
\label{corollary:projective_modules_over_the_twisted_group_algebra_of_e_are_free}
Let $U \leq \mathbf{G}_{n}$ be an open subgroup of the Morava stabilizer group which is a pro-$p$-group and which lies in the kernel of $\mathbf{G}_{n} \rightarrow \mathrm{Gal}(\mathbf{F}_{p^{n}} / \mathbf{F}_{p})$. Then any projective $E^{\ast} \llbracket U \rrbracket$ module is free in the graded sense.  
\end{corollary}

\begin{proof}
By a reduction to the ungraded case as in the proof of \cref{proposition:twisted_group_algebra_of_a_uniform_subgroup_is_finite_dim_and_notherian}, it is enough to show that any (ungraded) projective module over $E^{0} \llbracket U \rrbracket$ is free. As the latter is local by \cref{proposition:twisted_group_algebra_of_a_pro_p_group_is_local}, this follows from a theorem of Kaplansky that over a local ring any projective module is free \cite[{Theorem 2}]{kaplansky1958projective}.
\end{proof}

\section{Homotopy of \texorpdfstring{$K(n)$}{K(n)}-local finite complexes} 

In this section we construct a $K(n)$-local finite complex $V$ whose Adams-Novikov spectral sequence collapses on the second page and can be completely understood. The complex will be of the form 
\[
V \colonequals M \otimes E^{hU} 
\]
for a suitable generalized Moore spectrum 
\[
M=\thesphere/(p^{\alpha_0},v_1^{\alpha_1},\dots,v_{n-1}^{\alpha_{n-1}})
\]
and an open normal subgroup $U\leq\mathbf G_n$. In other words, it is a base-change of a generalized Moore spectrum to a finite Galois extension $E^{hU}$ of the $K(n)$-local sphere. The Adams-Novikov spectral sequence for $V$ takes the form
\[
    E_2^{s,t}=\Hrm^s(U,\pi_t(M\otimes E)) \Rightarrow \pi_{t-s}(V).
\]
In \cref{lemma:ANSS_E2}, we study the $E_2$-page of this spectral sequence and how it behaves when we shrink the subgroup $U$. In \cref{theorem:existence_of_a_convenient_generalized_moore_and_open_subgroup}, we prove that the spectral sequence collapses on the $E_2$-page for an appropriate choice of $M$ and $U$, and we show that $\pi_\ast V$ is an exterior algebra.

\begin{notation}
\label{notation:multi_indices_defining_generalized_moore_spectra}
In this section, let us write
\[
    \varepsilon:=\begin{cases}
2&\text{if $p=2$}\\
1&\text{if $p>2$.}
    \end{cases}
\]

If $(k, \alpha_{1}, \ldots, \alpha_{n-1}) \in \NN^{\times n}$ is a multi-index, we write 
\[
I(k, \alpha) \colonequals (p^{\varepsilon k}, v_{1}^{\alpha_{1}}, \ldots, v_{n-1}^{\alpha_{n-1}}) \subseteq E_{*}
\]
for the corresponding ideal given by setting $\alpha_{0} = \varepsilon k$. If it exists, we write 
\[
\thesphere / I(k, \alpha) \colonequals  \thesphere / (p^{\varepsilon k}, v_{1}^{\alpha_{1}}, \ldots, v_{n-1}^{\alpha_{n-1}})
\]
for the corresponding generalized Moore spectrum. Given an open subgroup $U\leq\mathbf G_n$, we write $E_r^{s,t}(U,k,\alpha)$ (not to be confused with Morava $E$-theory!) for the Adams-Novikov spectral sequence
\[
    E_r^{s,t}(U,k,\alpha)\Rightarrow \pi_{t-s}(E^{hU}/I(k,\alpha))
\]
whose $E_2$-page is given by
\[
    E_2^{s,t}=H^s(U,\pi_t(E/I(k,\alpha)))
\]
as previously discussed.
\end{notation}

\begin{remark}
It is likely that the use of $\varepsilon$, i.e. the use of powers of $p^2$ instead of $p$ at $p=2$, is unnecessary, but we use it to cite results from \cite{burklund2022multiplicative} directly.
\end{remark}

\begin{lemma}
\label{lemma:moore_multiplication}
    There exists an integer $k_0$ such that all of the following statements hold.
    \begin{enumerate}
        \item For any integer $k\geq k_0$, there exists an $\mathbf E_3$-ring structure on the Moore spectrum $\thesphere/p^{\varepsilon k}$ such that all of the projection maps
        \[
            \thesphere/p^{\varepsilon(k+1)}\to\thesphere/p^{\varepsilon k}
        \]
        are maps of $\mathbf E_3$-rings.
        \item For any integer $k\geq k_0$, there is a multi-index $\alpha\in\NN^{n-1}$ (depending on $k$) such that the generalized Moore spectrum
        \[
            \thesphere/I(k,\alpha)
        \]
        is an $\mathbf E_2$-$\thesphere/p^{\varepsilon k}$-algebra. Furthermore, for the same $\alpha$ and for all $k\geq \ell\geq k_0$, the generalized Moore spectrum
        \[
            \thesphere/I(\ell,\alpha)
        \]
        exists as an $\mathbf E_2$-$\thesphere/p^{\varepsilon \ell}$-algebra and the projection
        \[
            \thesphere /I(k,\alpha)\to\thesphere /I(\ell,\alpha)
        \]
        is a map of $\mathbf E_2$-$\thesphere /p^{\varepsilon k}$-algebras.
    \end{enumerate}
\end{lemma}
\begin{proof}
    Let us take $k_0=\max(4,n+2)$. Then, by \cite[Theorem 1.5]{burklund2022multiplicative}, there exists a tower of $\mathbf E_{\max(3,n+1)}$-rings
    \[
        \cdots\to\thesphere/p^{\varepsilon(k_0+1)}\to\thesphere/p^{\varepsilon(k_0)}
    \]
    implying (1). The first part of (2) follows from a repeated use of the same cited theorem, and it is essentially the content of \cite[Theorem 1.3]{burklund2022multiplicative}. The second part of (2) follows from the $\mathbf E_3$-structures of the projection maps in (1).
\end{proof}

\begin{lemma}
\label{lemma:small_enough_subgps_of_gn_act_trivially_on_homology_of_type_h_moore_spectra}
Let $\thesphere / I(k, \alpha)$ be a generalized Moore spectrum of type $n$. Then there exists an open subgroup $U \leq \mathbf{G}_{n}$ acting trivially on 
\[
E_{*}(\thesphere / I(k, \alpha)) \simeq E_{*}/I(k, \alpha). 
\]
\end{lemma}

\begin{proof}
Since $E_{*}/I(k, \alpha)$ is finitely generated as a ring, it is enough to show that there is an open subgroup stabilizing a finite set of generators. This is clear, as $E_{*}/I(k, \alpha)$ has a discrete topology and the action of $\mathbf{G}_{n}$ is continuous. 
\end{proof}

\begin{definition}
\label{definition:katrivial}
    Given a multi-index $(k,\alpha)\in\NN^n$ satisfying \cref{lemma:moore_multiplication}.(2), we say that an open subgroup $U\leq \mathbf G_n$ is \emph{$(k,\alpha)$-trivial} if
    \begin{enumerate}
        \item $U$ is $\varepsilon k$-uniform in the sense of \cref{definition:n_abelian_uniforn_group},
        \item $U$ acts trivially on $E_\ast/I(k,\alpha)$, and
        \item there exists a section of
        \[
            \pi_\ast(E^{hU}/I(k,\alpha))\to \pi_\ast (E/I(k,\alpha))=E_\ast/I(k,\alpha)
        \]
        as a map of rings.
    \end{enumerate}
\end{definition}

\begin{remark}
As a consequence of \cref{lemma:subgroup_of_pth_powers_of_an_n_abelian_gp_is_nplusone_abelian}, if $U$ is $(k, \alpha)$-trivial, then so are all of the subgroups $U_{t}$ in its lower $p$-series. Also, if $U$ is $(k,\alpha)$-trivial, then it is $(k',\alpha)$-trivial for any $k_0\leq k'\leq k$. Indeed, given a section $E_\ast/I(k,\alpha)\to \pi_\ast(E^{hU}/I(k,\alpha))$, we obtain a section
\[
    E_\ast/I(k',\alpha)=(E_\ast/I(k,\alpha))/p^{\varepsilon k'}\to \pi_\ast(E^{hU}/I(k,\alpha))/p^{\varepsilon k'}\to \pi_\ast (E^{hU}/I(k',\alpha))
\]
where the last map comes from the fact that $p^{\varepsilon k'}=0$ in the last group as the relevant ring spectrum is a $\thesphere/p^{\varepsilon k'}$-algebra. 
\end{remark}

\begin{lemma}
\label{lemma:ANSS_E2}
    For any $k$ and $\alpha$ as in \cref{lemma:moore_multiplication}.(2), there exists a $(k,\alpha)$-trivial open normal subgroup $U \lhd \mathbf G_n$. Also, being $(k,\alpha)$-trivial implies all of the following, where we write $I=I(k,\alpha)$:
    \begin{enumerate}
        \item For any $m\geq1$, there is an isomorphism of $E_\ast/I$-algebras
        \[
            \Hrm^\ast(U_m,E_\ast/I) \simeq\Lambda_{E_\ast/I}(x_1,\dots,x_{n^2})
        \]
        for any basis $\{x_i\}\subseteq \Hrm^1(U_m,\ZZ/p^{\varepsilon k})$, whose image in $\Hrm^1(U_m,E_\ast/I)$ is denoted in the same way, where $U_{m} = \{ u^{p^{m-1}} \ | \ u \in U\}$ is as in  \cref{notation:lower_p_series}.
        \item The map
        \[
            \Hrm^i(U_m,E_\ast/I)\to \Hrm^i(U_{m+\varepsilon},E_\ast/I)
        \]
        induced by the inclusion $U_{m+\varepsilon}\to U_m$ is multiplication by $p^{\varepsilon i}$ for a suitable choice of basis for both sides. Note that both sides are free $E_\ast/I$-modules.
    \end{enumerate}
    The above statements hold for any $k_0\leq k'<k$ with the same $U$ and $\alpha$, compatibly with the natural maps $H^\ast(U_m,E_\ast/I(k,\alpha))\to H^\ast(U_m, E_\ast/I(k',\alpha))$. 
\end{lemma}
\begin{proof}
Since there exists an open pro-$p$ subgroup of $\mathbf G_n$, we have
\[
\varinjlim \pi_{*}(E^{hU}/I) \simeq E_{*}/I \simeq \mathrm{W}(\mathbf{F}_{p^{n}})[v_{1}, \ldots, v_{n-1}][u^{\pm 1}]/(p^{\varepsilon k}, v_{1}^{\alpha_{1}}, \ldots, v_{n-1}^{\alpha_{n-1}})
\]
where the colimit is taken over open normal pro-$p$ subgroups $U\lhd\mathbf G_n$. The right hand side is a finitely presented ring, hence the map from $\pi_{*}(E^{hU}/I)$ admits a section for all sufficiently small subgroups $U$. By \cref{lemma:small_enough_subgps_of_gn_act_trivially_on_homology_of_type_h_moore_spectra}, we can make $U$ smaller so that it acts trivially on $E_\ast/I$.

Since $U$ is a $p$-adic analytic pro-$p$ group, the lower $p$-series subgroups $U_{k}$ are uniform for sufficiently large $k$; see \cite[Proposition 3.9 and Theorem 4.5]{dixon2003analytic}. Further replacing $U$ by $U_{k+m}$ we can assume that it is $m$-uniform for any chosen $m$ by \cref{lemma:subgroup_of_pth_powers_of_an_n_abelian_gp_is_nplusone_abelian}. This proves the first statement. 

Then, part (1) of the statement follows from \cref{proposition:zpn_cohomology_of_n_abelian_uniform_group_is_exterior} since $E_\ast/I$ is finite free over $\ZZ/p^{\varepsilon k}$ in each (topological) degree.

For (2), it is enough to consider $i=1$, since higher cohomology groups are multiplicatively generated by $\Hrm^1$. Recall from the proof of \cref{proposition:zpn_cohomology_of_n_abelian_uniform_group_is_exterior} that there are natural isomorphisms
    \[
        \Hrm^1(U_m,E_\ast/I) \simeq \mathrm{Hom}(U_m,E_\ast/I)\simeq \mathrm{Hom}(U_m/U_{m+\varepsilon k},E_\ast/I) =\mathrm{Hom}((\ZZ/p^{\varepsilon k})^d,E_\ast/I)
    \]
    for any $m$. Under these isomorphisms, the map in question corresponds to
    \[
        \mathrm{Hom}(U_m/U_{m+\varepsilon k},E_\ast/I)\to\mathrm{Hom}(U_{m+\varepsilon}/U_{m+\varepsilon k},E_\ast/I)\to\mathrm{Hom}(U_{m+\varepsilon}/U_{m+\varepsilon(1+k)},E_\ast/I)
    \]
    induced by the obvious inclusion and projection. Since $U_{m+\varepsilon}$ is the subgroup of $p^{\varepsilon}$-th powers in $U_m$, the conclusion easily follows.

    The last assertion follows from the fact that $U$ is $(k',\alpha)$-trivial.
\end{proof}

Next, we turn our attention to the Adams-Novikov spectral sequence
\[
    E_2^{s, \ast}(U,k,\alpha)=H^s(U,E_\ast/I(k,\alpha))\Rightarrow \pi_{\ast-s} (E^{hU}/I(k,\alpha)).
\]
Assuming $k,\alpha,U$ are chosen to satisfy \cref{lemma:ANSS_E2}, the spectral sequence has commutative multiplication and we understand the $E_2$-page as a ring.

\begin{lemma}
\label{lemma:ANSS_degeneration}
Let $(k,\alpha)\in\NN^n$ be a multi-index satisfying (2) of  \cref{lemma:moore_multiplication} and $U$ a $(k,\alpha)$-trivial open subgroup of $\mathbf G_n$. Also, let $q\geq1$ be an integer and assume that $k/(q+1)$ is an integer greater than or equal to $k_0$ of \cref{lemma:moore_multiplication}. If the spectral sequence $E_r(U,k,\alpha)$ has no $d_i$-differentials for $2\leq i\leq q$, then there exists a subgroup $V \leq U$ in the lower $p$-series of $U$ (in particular, characteristic) which is $(k,\alpha)$-trivial and such that the spectral sequence $E_r(V,k/(q+1),\alpha)$ has no $d_i$-differentials for $2\leq i\leq q+1$.
\end{lemma}  

\begin{proof}
We claim that $V=U_{1+\varepsilon t}$ for $t:=k/(q+1)$ has the required properties. To show that $E_r(V,t,\alpha)$ has no $d_i$-differentials for $2\leq i\leq q+1$, it is enough to show that the classes in
    \[
        E_2^{s, \ast}(V,t,\alpha)=\Hrm^s(V,E_\ast/I(t,\alpha))\qquad s=0,1
    \]
    are $d_i$-cycles since the $E_2$-page is multiplicatively generated by these two lines.

    All classes on the zeroth line
    \[
        E_2^{0,\ast}(V,t,\alpha)=\mathrm{H}^0(V,E_\ast/I(t,\alpha))=E_\ast/I(t,\alpha)
    \]
    are permanent cycles because they must detect the image of the section
    \[
        E_\ast/I(t,\alpha)\to \pi_\ast (E^{hV}/I(t,\alpha)).
    \]

    Next, we study classes on the first line, i.e. in $E_2^{1, \ast}(V,t,\alpha)$. For $d_i$-differentials for $2\leq i\leq q$, we prove the following statement using induction on $j$: for any $2\leq j\leq q+1$, every class of $E_2^{1, \ast}(V,k-(j-2)t,\alpha)$ is a $d_i$-cycle for $2\leq i\leq j-1$.

    For $j=2$, there is nothing to prove. Now, assume that the statement holds for some $j\leq q$. Then, considering the map of spectral sequences
    \[
        E_r(V,k-(j-2)t,\alpha)\to E_r(V,k-(j-1)t,\alpha),
    \]
    every class of $E_2^{1, \ast}(V,k-(j-1)t,\alpha)$ is a $d_i$-cycle for $2\leq i<j$.
    
    To show that these classes are $d_j$-cycles, consider maps of spectral sequences
    \[
        E_r(U,k-(j-2)t,\alpha)\to E_r(V,k-(j-2)t,\alpha)\to E_r(V,k-(j-1)t,\alpha).
    \]
    Let $x\in E_2^{1, \ast}(V,k-(j-1)t,\alpha)$ be any class. This class lifts along the second map to a class in $E_2^{1, \ast}(V,k-(j-2)t,\alpha)$, which we also denote by $x$, since that map is a reduction modulo $p^{\varepsilon(k-(j-1)t)}$. By \cref{lemma:ANSS_E2}, on $E_2^{1, \ast}$, the first map above is multiplication by $p^{\varepsilon t}$. In particular,
    \[
        p^{\varepsilon t}x\in E_2^{1, \ast}(V,k-(j-2)t,\alpha)
    \]
    is in the image of the first map, and therefore is a $d_j$-cycle by our assumption of the lemma since $j\leq q$. By the induction hypothesis and multiplicativity, we have
    \[
        E_2^{\ast, \ast}(V,k-(j-2)t,\alpha)=E_j^{\ast, \ast}(V,k-(j-2)t,\alpha)
    \]
    and in particular, the right hand side is free over $\ZZ/p^{\varepsilon(k-(j-2)t)}$. Therefore, $p^{\varepsilon t}x$ being a $d_j$-cycle implies that
    \[
        d_j(x)\in p^{\varepsilon(k-(j-1)t)}\cdot E_j^{j+1,\ast}(V,k-(j-2)t,\alpha)
    \]
    which becomes zero after passing to $E_j^{\ast, \ast}(V,k-(j-1)t,\alpha)$.

    This finishes the induction and considering the case $j=q+1$, we have proved that every class of $E_2^{1, \ast}(V,2t,\alpha)$ is a $d_i$-cycle for $2\leq i\leq q$. Finally, consider the maps of spectral sequences
    \[
        E_r(U,2t,\alpha)\to E_r(V,2t,\alpha)\to E_r(V,t,\alpha).
    \]
    On the $E_2$-page, the first map is multiplication by $p^{\varepsilon ts}$ on $E_2^{s, \ast}$, which in particular is zero if $s\geq 2$.
    
    The rest of the argument is similar to the previous cases. Given a class $x\in E_2^{1, \ast}(V,t,\alpha)$, we can lift it to a class $x\in E_2^{1, \ast}(V,2t,\alpha)$, and by \cref{lemma:ANSS_E2}, the class $p^{\varepsilon t}x$ is in the image of the first map. Therefore, it must be a $d_{q+1}$-cycle since the first map is the zero map above the first line. Since $E_{q+1}(V,2t,\alpha)$ is free over $\ZZ/p^{2\varepsilon t}$, it implies that $d_{q+1}(x)$ is zero mod $p^{\varepsilon t}$.
\end{proof}

\begin{theorem}
\label{theorem:existence_of_a_convenient_generalized_moore_and_open_subgroup}
For any integer $k\geq k_0$ (\cref{lemma:moore_multiplication}), there exists a multi-index $\alpha\in\NN^{n-1}$ and a $(k,\alpha)$-trivial open normal subgroup $U \lhd \mathbf G_n$ such that the following statements are satisfied, where we write $I=I(k,\alpha)$ for brevity.
    \begin{enumerate}
        \item The Adams-Novikov spectral sequence
        \[
            E_2(U,k,\alpha)=H^\ast(U,E_\ast/I)\Rightarrow \pi_\ast(E^{hU}/I)
        \]
        collapses on the second page.
        \item For any collection of lifts $\overline{x}_{i} \in \pi_{-1}(E^{hU}/I)$ of $x_{i}$ there exists an isomorphism
\[
\pi_{\ast}(E^{hU}/I) \simeq \Lambda_{E_{*}/I}(\overline{x}_{1}, \ldots, \overline{x}_{n^2}) 
\]
of filtered algebras, where we equip the left hand side with the Adams-Novikov filtration and the right hand side with the decreasing filtration induced from the grading where $E_{*}/I$ is of degree zero and the $\overline{x}_{i}$'s are of degree one. 
    \end{enumerate}
\end{theorem}

\begin{proof}
Let $K=(n^2)! \cdot k$ and choose $\alpha$ that satisfies \cref{lemma:moore_multiplication} and a $(K,\alpha)$-trivial open normal subgroup $V \lhd \mathbf{G}_{n}$. Repeatedly using \cref{lemma:ANSS_degeneration} for $q=1,\dots,n^2-1$, we obtain a $(k,\alpha)$-trivial $U \lhd V$ such that the Adams-Novikov spectral sequence
\begin{equation}
\label{equation:adams_novikov_collapsing_for_a_suitable_k_alpha_trivial_subgroup}
E_2^{s, \ast}=\Hrm^s(U,E_\ast/I)\Rightarrow \pi_\ast(E^{hU}/I)
\end{equation}
has no $d_i$-differentials for $i\leq n^2$. This implies that this spectral sequence collapses on the $E_2$-page since there is a horizontal vanishing line at $s=n^2$ on the $E_2$-page. Note that $U$ is also normal in $\mathbf G_n$ since it is a characteristic subgroup of $V$ by construction.

Now suppose we are given lifts $\overline{x}_{i}$ of $x_{i}$. If $p > 2$ or $n = 1$, these lifts satisfy $\overline{x_{i}}^{2} = 0$ by anticommutativity or a vanishing line argument. Together with the section $E_\ast/I\to \pi_{*}(E^{hU}/I)$ which exists by $(k, \alpha)$-triviality these lifts define a map of filtered algebras 
\[
\Lambda_{E_{*}/I}(\overline{x}_{1}, \ldots, \overline{x}_{n^{2}}) \rightarrow \pi_{*}(E^{hU}/I). 
\]
This is an isomorphism on associated graded pieces by construction and hence an isomorphism. 

The argument from the previous paragraph also works if $p=2$ and $n > 1$ if we can construct lifts of $x_{i}$ which satisfy $\overline{x}_{i}^{2} = 0$. We explain how to do so now. First, we rerun the argument for (1) with $k$ replaced by $\ell=k+n^2-2$. Then, we obtain a multi-index $\alpha$ and a $(\ell,\alpha)$-trivial subgroup $U\leq\mathbf G_n$ such that the $E_\infty$-page of the Adams-Novikov spectral sequence is given by
\[
    E_\infty^{\ast, \ast}(U,\ell,\alpha)=\Lambda_{E_\ast/I(\ell,\alpha)}(x_1,\dots,x_{n^2})\Rightarrow\pi_\ast(E^{hU}/I(\ell,\alpha)).
\]
Using the comparison map $E^{hU}/I(\ell,\alpha)\to E^{hU}/I(k,\alpha)$, it is easy to see that part (1) for $(U,\ell,\alpha)$ implies that it also holds for $(U,k,\alpha)$.

Given any lifts $\overline{x_i}\in\pi_{-1}(E^{hU}/I(k,\alpha))$ of the generators $x_i\in E_\infty^{\ast,\ast}(U,k,\alpha)$, we can further lift them to $\overline{x_i}\in\pi_{-1}(E^{hU}/I(\ell,\alpha))$, since we can see from the spectral sequences that $\pi_\ast(E^{hU}/I(\ell,\alpha))\to\pi_\ast(E^{hU}/I(k,\alpha))$ is surjective.

By induction on $j$, we shall prove that for all $i$, the class $\overline{x_i}^2$ has Adams-Novikov filtration $\geq 3+j$ in $\pi_\ast(E^{hU}/I(\ell-j,\alpha))$. Then, the proof is done if we consider $j=n^2-2$, by the vanishing line at $s=n^2$.

The above is true for $j=0$ since $x_i^2=0$ on the $E_\infty$-page. Suppose it is true for some $j$. If $\overline{x_i}^2=0$ in $\pi_\ast(E^{hU}/I(\ell-j,\alpha))$, then there is nothing to prove. Otherwise, let $y$ be a class in $E_\infty^{\ast, \ast}(U,\ell-j,\alpha)$ that detects $\overline{x_i}^2$. By anticommutativity, $\overline{x_i}^2$ is $2$-torsion. Therefore, $y$ is $2$-torsion, and it must map to zero under
\[
    E_\infty^{\ast, \ast}(U,\ell-j,\alpha)\to E_\infty^{\ast, \ast}(U,\ell-j-1,\alpha)
\]
since it is a mod $4^{\ell-j-1}$ reduction of a free $\ZZ/4^{\ell-j}$ module. This implies that $\overline{x_i}^2$ has higher Adams-Novikov filtration in $\pi_\ast(E^{hU}/I(\ell-j-1,\alpha))$ than in $\pi_\ast(E^{hU}/I(\ell-j,\alpha))$.
\end{proof}

\section{Characterization of locally fp-spectra} 

In this section, we prove the main result of this note, \cref{theorem:classification_of_local_fp_types_in_terms_of_cohomology_and_morava_e_theory}. Before we proceed, we recall that we call a $K(n)$-local spectrum $X$ \emph{locally fp} if $\pi_{*}(V \otimes X)$ is degreewise finite for any finite complex $V$ of type $n$. Our terminology is motivated by the following easy observation: 

\begin{lemma}
\label{lemma:spectra_of_fp_type_n_which_satisfy_the_telescope_localize_to_local_fps}
Let $X$ be a spectrum of fp-type $n$ in the sense of Mahowald and Rezk and suppose that $L_{T(n)} X \simeq L_{K(n)}X$. Then $L_{K(n)} X$ is locally fp. 
\end{lemma}

\begin{proof}
Let $V$ be a finite complex of type $n$ and choose a $v_{n}$-self map $v \colon V \rightarrow V$. By assumption, 
\[
\pi_{*} (V/v \otimes X) 
\]
is finite, and hence so is $\pi_{*}(V \otimes X)/v$. It follows that we can choose finitely many classes $x_{1}, \ldots, x_{m} \in \pi_{*}(V \otimes X)$ whose $v$-power multiples generate $\pi_{*}(V \otimes X)$ as an abelian group. Since $n > 0$, these classes are torsion, so that 
\[
v^{-1} \pi_{*}(V \otimes X) \simeq \pi_{*}(V[v^{-1}] \otimes X) \simeq \pi_{*}(L_{T(n)}(V \otimes X)) \simeq \pi_{*}(V \otimes L_{K(n)}X)
\]
is degreewise finite, as needed. 
\end{proof}

\begin{notation}
We write $- \widehat{\otimes} - \colonequals L_{K(n)}(- \otimes -)$ for the $K(n)$-local tensor product. If $X$ is $K(n)$-local, we write $D(X) \colonequals F(X, \thesphere_{K(n)})$ for its monoidal dual. 
\end{notation}

\begin{theorem}
\label{theorem:classification_of_local_fp_types_in_terms_of_cohomology_and_morava_e_theory}
Let $X$ be a $K(n)$-local spectrum. Then the following are equivalent: 
\begin{enumerate}
    \item $X$ is locally fp; that is, $\pi_{*}(V \otimes X)$ is degreewise finite for any type $n$ complex $V$, 
    \item $E^{*}X$ is finitely generated as an $E^{*}E$-module, 
    \item $X$ belongs to the thick subcategory of spectra generated by $E$, 
    \item $X$ is reflexive; that is, the canonical map $X \to D^{2}(X)$ into its $K(n)$-local monoidal double dual is an equivalence.
\end{enumerate}
\end{theorem}

The proof of \cref{theorem:classification_of_local_fp_types_in_terms_of_cohomology_and_morava_e_theory} will take the remainder of this section. As we explained in the introduction, of the four implications that make the proof, it is $(1 \Rightarrow 2)$ that is particularly delicate. To make our arguments easier to process, we first deal with the other three implications. 

\begin{lemma}
\label{lemma:fixed_points_for_an_open_normal_generate_contain_the_sphere_in_thick_ideal}
Let $U \leq \mathbf{G}_{n}$ be an open subgroup. Then $\thesphere_{K(n)}$ belongs to the thick tensor-ideal of dualizable $K(n)$-local spectra generated by $E^{hU}$. 
\end{lemma}

\begin{proof}
By \cite[{Corollary 2.7, Proposition 4.6}]{barthel2022conjectures}, it is enough to show that $E = E \widehat{\otimes} \thesphere_{K(n)}$ belongs to the thick ideal of compact $E$-modules generated by $E \widehat{\otimes} E^{hU}$. Since the homotopy of the latter
\[
E_{*}^{\wedge}(E^{hU}) \simeq \map_{\mathrm{cts}}(\mathbf{G}_{n}/U, E_{*})
\]
is a free $E_{*}$-module, $E \widehat{\otimes} E^{hU}$ is a free $E$-module, ending the argument. 
\end{proof}

\begin{lemma}
\label{lemma:being_in_thick_subcategory_of_e_up_to_finite_galois_extension_means_youre_in_thick_subcat_of_e}
Let $X$ be a $K(n)$-local spectrum such that $E^{hU} \widehat{\otimes} X$ belongs to the thick subcategory of $K(n)$-local spectra generated by $E$ for some open $U \leq \mathbf{G}_{n}$. Then $X$ also belongs to this thick subcategory. 
\end{lemma}

\begin{proof}
Let $\ccat \subseteq \spectra_{K(n)}^{\mathrm{dual}}$ denote the subcategory of those dualizable $K(n)$-local spectra $A$ such that $A \widehat{\otimes} X$ belongs to the thick subcategory generated by $E$. By \cref{lemma:fixed_points_for_an_open_normal_generate_contain_the_sphere_in_thick_ideal}, it is enough to verify that $\ccat$ is a tensor-ideal. By taking $M = A \widehat{\otimes} X$, it is sufficient to check that if $M$ belongs to the thick subcategory generated by $E$, then so does $B \widehat{\otimes} M$ for an arbitrary dualizable $B$. 

If we fix $B$, the collection of $M$ for which this holds is clearly thick, so it is enough to check the case $M = E$. Since $\pi_{*}(B \widehat{\otimes} E) = E_{*}^{\wedge}(B)$ is a finitely generated $E_{*}$-module over the regular ring $E_{*}$ when $B$ is dualizable \cite[{Theorem 8.6}]{hovey1999morava}, it admits a finite resolution by finite projective $E_{*}$-modules. Lifting this resolution to $E$-modules shows that $B \widehat{\otimes} E$ is in the thick subcategory generated by $E$ as needed. 
\end{proof}

\begin{proof}[{Proof of \cref{theorem:classification_of_local_fp_types_in_terms_of_cohomology_and_morava_e_theory}, implications $(2 \Rightarrow 3)$, $(3 \Rightarrow 4)$ and $(4 \Rightarrow 1)$:}] 
We first show $(2 \Rightarrow 3)$. If $E^{\ast}X$ is finitely generated, a choice of a finite set of generators determines a homotopy class of a map of spectra
\[
X \rightarrow F
\]
which is surjective on $E$-cohomology and such that $F$ is a finite sum of shifts of $E$. Let $C \colonequals \cof(X\to F)$. Since $E^{\ast}E$ is noetherian by \cref{corollary:twisted_group_algebra_of_gn_is_noether}
\[
E^{*}C \simeq \mathrm{ker}(E^{*}F \rightarrow E^{*}X) 
\]
is also finitely generated. By construction, $C$ belongs to the thick subcategory generated by $E$ if and only if $X$ does. 

Let $U$ be a uniform open subgroup of $\mathrm{ker}(\mathbf{G}_{n} \rightarrow \mathrm{Gal}(\mathbf{F}_{p^{n}} / \mathbf{F}_{p}))$. Since $E^{*}E$ is free as an $E^{*} \llbracket U \rrbracket$-module, so is $E^{*}F$. It follows that projective dimension of the cohomology of $C$ satisfies 
\[
\mathrm{pd}_{E^{*} \llbracket U \rrbracket}(E^{*}C) = \mathrm{max}(0, \mathrm{pd}_{E^{*} \llbracket U \rrbracket}(E^{*}X)-1).
\]
Since $E^{*} \llbracket U \rrbracket$ is of finite global dimension $n^{2}+n$ by \cref{proposition:twisted_group_algebra_of_a_uniform_subgroup_is_finite_dim_and_notherian}, repeating the above operation $n^{2}+n$ times we can assume that $E^{*}X$ is projective as an $E^{*} \llbracket U \rrbracket$-module, and hence free by \cref{corollary:projective_modules_over_the_twisted_group_algebra_of_e_are_free}. 

Using the K\"{u}nneth formula, we have
\[
E^{*}(E^{hU} \widehat{\otimes} X) \simeq E^{*}E^{hU} \otimes_{E_{*}} E^{*}X \simeq E^{*}[\mathbf{G}_{n}/U] \otimes_{E_{*}} E^{*}X,
\]
where the completed tensor product is the $K(n)$-local one. As an $E^{*}E$-module, we can identify the right hand side with the extension of scalars
\[
E^{*}(E^{hU} \widehat{\otimes} X) \simeq E^{*}E \otimes_{E^\ast\llbracket U\rrbracket } E^{*}X. 
\]
Hence the left hand side is a free $E^{*}E$-module of finite rank as $E^{*}X$ was assumed to be free over $E^\ast\llbracket U\rrbracket$. A choice of a basis determines a homotopy class of a map $E^{hU} \widehat{\otimes} X \rightarrow F$ into a free $E$-module. By construction, this map is an $E^{*}(-)$-isomorphism and hence an equivalence as both sides are $K(n)$-local \cite[Corollary 7.6]{hovey1999morava}. It follows that $E^{hU} \widehat{\otimes} X$ is in the thick subcategory generated by $E$, and hence so is $X$ itself by \cref{lemma:being_in_thick_subcategory_of_e_up_to_finite_galois_extension_means_youre_in_thick_subcat_of_e}. 

The implication $(3 \Rightarrow 4)$ follows from the fact that reflexive spectra form a thick subcategory and that $E$ is reflexive by a result of Strickland \cite[{Proposition 16}]{strickland2000gross}. 

To show $(4 \Rightarrow 1)$, we make use of Brown-Comenetz duality. We write 
\[
\widehat{I} \colonequals F(\mathrm{fib}(L_{n} \thesphere \rightarrow L_{n-1} \thesphere), I_{\mathbf{Q}/\mathbf{Z}})
\]
for the $K(n)$-local Brown-Comenetz dual of the sphere, which is an invertible $K(n)$-local spectrum by \cite[{Theorem 10.2}]{hovey1999morava}, which implies that $F(X, \widehat{I}) \simeq D(X) \widehat{\otimes} \widehat{I}$ for any $K(n)$-local spectrum $X$. It follows that if $X$ is reflexive, then the Brown-Comenetz double dual map
\[
X \rightarrow F(F(X, \widehat{I}), \widehat{I})
\]
is also an equivalence. Tensoring both sides with a finite complex of type $n$, we obtain an equivalence 
\[
V \otimes X \rightarrow F(F(V \otimes X, \widehat{I}), \widehat{I}). 
\]
Since $V \otimes X$ is $L_{n}$-local and $L_{n-1}$-acyclic, this map can be identified with the integral Brown-Comenetz dual map 
\[
V \otimes X \rightarrow F(F(V \otimes X, I_{\mathbf{Q}/\mathbf{Z}}), I_{\mathbf{Q}/\mathbf{Z}}), 
\]
which at the level of homotopy groups we can identify with the Pontryagin double dual 
\[
\pi_{k}(V \otimes X) \rightarrow \Hom_{\abeliangroups}(\Hom_{\abeliangroups}(\pi_{k}(V \otimes X), \mathbf{Q}/\mathbf{Z}), \mathbf{Q}/\mathbf{Z}). 
\]
Suppose, by contradiction, that $X$ is not locally fp, so that there exists a $k \in \ZZ$ such that $\pi_{k}(V \otimes X)$ is infinite. Since $n > 0$, $\pi_{k}(V \otimes X)$ is a $\ZZ/p^{a}$-module for some $a$, so that $\mathbf{F}_{p} \otimes_{\ZZ} \pi_{k}(V \otimes X)$ is also infinite. Then 
\[
\mathbf{F}_{p} \otimes_{\ZZ} \Hom(\Hom(\pi_{k}(V \otimes X), \mathbf{Q}/\mathbf{Z}), \mathbf{Q}/\mathbf{Z}) \simeq \Hom(\Hom(\mathbf{F}_{p} \otimes_{\ZZ} \pi_{k}(V \otimes X), \mathbf{Q}/\mathbf{Z}), \mathbf{Q}/\mathbf{Z})
\]
is an $\mathbf{F}_{p}$-vector space of strictly larger cardinality, which contradicts the Pontryagin double dual map being an isomorphism. This ends the proof. 
\end{proof}

The proof of the implication $(1 \Rightarrow 2)$ is a bit longer, and will make use of finite ring spectra constructed in the previous section. We first set up the notation. 

\begin{notation}
\label{notation:convenient_type_h_complex}
In \cref{theorem:existence_of_a_convenient_generalized_moore_and_open_subgroup}, we had shown the existence of a particularly well-behaved pair of a generalized Moore spectrum $M$ of type $n$ and an open subgroup $U \lhd \mathbf{G}_{n}$ of the Morava stabilizer group.  We write 
\[
V = M \otimes E^{hU}, 
\]
which is compact and hence a retract of the $K(n)$-localization of a finite complex of type $n$. This is an $\mathbf{E}_{2}$-$E^{hU}$-algebra and we have an isomorphism 
\[
\pi_{*}(V) \simeq \Lambda_{E_{*}/I}(x_{1}, \ldots, x_{d})
\]
for some $x_{i} \in \pi_{-1}(V)$ and $d = n^{2}$, and the Adams-Novikov filtration coincides with the filtration induced by the grading with $x_{i}$ of degree one. Thus, $V_{*} = \pi_{*}(V)$ is a $2$-periodic Artin local Dirac ring\footnote{Here, following \cite{hesselholt2023dirac, hesselholt2024dirac}, we refer to graded-commutative rings satisfying the Koszul sign rule as \emph{Dirac rings.} In a Dirac ring, left and right graded ideals coincide by \cite[{Lemma 2.5}]{hesselholt2023dirac}. A Dirac ring is \emph{local} if it has a single maximal graded ideal and it is \emph{Artin} if graded ideals satisfy the descending chain condition.} with residue field 
\[
K_{*} \colonequals E_{*}/\mfrak \simeq \mathbf F_{p^n}[u^{\pm 1}],
\]
where $\mfrak = (v_{0}, \ldots, v_{n-1}) \subseteq E_{*}$ is the maximal ideal. 
\end{notation}

\begin{lemma}
\label{lemma:characterization_of_fp_modules_over_homotopy_of_generalized_moore}
Let $N$ be a $V_{*}$-module. Then the following are equivalent: 
\begin{enumerate} 
    \item $N$ can be obtained as an iterated extension of shifts $K_{*}$ and $K_{*}[1]$ of the residue field, 
    \item $N$ is finitely presented, 
    \item $N$ is degreewise finite. 
\end{enumerate}
\end{lemma}

\begin{proof}
Both of the forward implications are clear, so we only show $(3) \Rightarrow (1)$. We write 
\[
\mathfrak{n} \colonequals (\mfrak, x_{1}, \ldots, x_{n^2}) \subseteq V_{*}
\]
for the unique maximal ideal, which is nilpotent as $V_{*}$ is Artin. Suppose that $N$ is degreewise finite. We have a surjection
\[
N \rightarrow N/\mathfrak{n} N \simeq K_{*} \otimes_{V_{*}} N.  
\]
If the target is zero, then $\mathfrak{n} N = N$. We deduce  by induction that $\mathfrak{n}^{k} N = N$ for all $k \geq 0$, which by nilpotence shows that $N = 0$. 

If $K_{*} \otimes_{V_{*}} N$ is non-zero, then it is a degreewise finite vector space over the $2$-periodic graded field $K_{*}$, so it is a finite sum of $K_{*}$ and $K_{*}[1]$. Replacing $N$ by the kernel of the reduction map we reduce the result to the case of a module with strictly smaller sum of cardinalities in degree zero and one. Descending induction on cardinality finishes the argument. 
\end{proof}

\begin{lemma}
\label{lemma:ext_algebra_noether_and_finite_presentation_over_ext_algebra}
The $\Ext$-algebra 
\begin{equation}
\label{equation:ext_algebra_of_homotopy_of_the_convenient_complex}
\Ext^{\ast, \ast}_{V_{*}}(E_{*}/I, E_{*}/I)
\end{equation}
is noetherian and concentrated in even total degree. Moreover, for any finitely presented $V_{*}$-module $M$, 
\[
\Ext^{\ast, \ast}_{V_{*}}(M, E_{*}/I)
\]
is a finitely presented module over (\ref{equation:ext_algebra_of_homotopy_of_the_convenient_complex}). 
\end{lemma}

\begin{proof}
The first part is clear, since $V_{*}$ is exterior over $E_{*}/I$, so that 
\[
\Ext^{\ast, \ast}_{V_{*}}(E_{*}/I, E_{*}/I) \simeq (E_{*}/I)[y_{1}, \ldots, y_{d}] 
\]
is a polynomial $E_{*}/I$-algebra on classes $y_{i}$ dual to $x_{i}$. For the second part, using the first characterization of \cref{lemma:characterization_of_fp_modules_over_homotopy_of_generalized_moore}, it is enough to show the result for $M = K_{*}$. Since $E_{*}/I$ is a self-injective ring as it is a complete intersection Artin local ring, we have
\[
\Ext^{\ast, \ast}_{E_{*}/I}(K_{*}, E_{*}/I) \simeq K_{*}, 
\]
concentrated in cohomological degree zero. It follows that the change-of-rings spectral sequence
\[
    E_2=\Ext^{\ast,\ast}_{E_\ast/I}(K_\ast, \Ext_{V_\ast}^{\ast,\ast}(E_\ast/I,E_\ast/I)) \Rightarrow \Ext^{\ast,\ast}_{V_\ast}(K_\ast,E_\ast/I)
\]
collapses on the second page, which is isomorphic to
\[
     \Ext^{\ast,\ast}_{E_\ast/I}(K_\ast, (E_{*}/I)[y_{1}, \ldots, y_{d}]) \simeq K_\ast[y_1,\dots,y_d]
\]
concentrated in outer cohomological degree zero. Therefore,
\[
\Ext^{\ast, \ast}_{V_\ast}(K_{*}, E_{*}/I) \simeq K_{*}[y_{1}, \ldots, y_{d}] 
\]
is a cyclic module over the needed $\Ext$-algebra. 
\end{proof}

\begin{lemma}
\label{lemma:collapse_of_uct_for_emodi_on_second_page}
Consider $E \otimes M$, the tensor product of Morava $E$-theory and the generalized Moore spectrum of \cref{notation:convenient_type_h_complex}, so that $\pi_{*}(E \otimes M) \simeq E_{*}/I$. We consider this tensor product as a left $V$-module using the map of $\mathbf{E}_{2}$-algebras 
\[
V = E^{hU} \otimes M \rightarrow E \otimes M 
\]
induced by the canonical map $E^{hU} \rightarrow E$. Then the associated universal coefficient spectral sequence 
\begin{equation}
\label{equation:uct_sseq_for_e_tensor_generalized_moore}
\Ext_{V_{*}}^{\ast}(E_{*}/I, E_{*}/I) \Rightarrow \pi_{*} \map_{\Mod_{V}}(E \otimes M, E \otimes M) 
\end{equation}
collapses on the second page and converges strongly.
\end{lemma}

\begin{proof}
By \cref{lemma:ext_algebra_noether_and_finite_presentation_over_ext_algebra}, the second page of the spectral sequence is concentrated in even total degree. Since the differentials are all of total odd degree, the collapse follows. 

The universal coefficient spectral sequence is conditionally convergent \cite[Theorem IV.4.1]{elmendorf2007rings}, and since it collapses on the second page, it is actually strongly convergent by a result of Boardman \cite[Remark following Theorem 7.1]{boardman1999conditionally}. 
\end{proof}

\begin{remark}
\label{remark:self_maps_of_emodi_as_an_f_module}
The abutment of the spectral sequence of \cref{lemma:collapse_of_uct_for_emodi_on_second_page} can be described explicitly. We have 
\[
\pi_{-*} \map_{\Mod_{E^{hU}}}(E, E) \simeq E^{*} \llbracket U \rrbracket,
\]
by \cref{lemma:ehu_linear_endo_of_e_is_a_twisted_group_algebra}, which is an $E_{*}$-module for which the sequence $v_{0}^{\alpha_{0}}, \ldots, v_{n-1}^{\alpha_{n-1}}$ of generators of $I$ is regular. It follows that the map of rings 
\[
\pi_{*} \map_{\Mod_{E^{hU}}}(E, E) \rightarrow \pi_{*} \map_{\Mod_{V}}(E \otimes M, E \otimes M) 
\]
induced by extension of scalars along $E^{hU} \rightarrow V$ is a surjection and induces an isomorphism 
\[
\pi_{-*} \map_{\Mod_{V}}(E \otimes M, E \otimes M) \simeq (E^{*}/I) \llbracket U \rrbracket. 
\]
Comparing their associated graded, it is plausible that the filtration induced by the spectral sequence of \cref{lemma:collapse_of_uct_for_emodi_on_second_page} coincides with the one induced by the lower $p$-series of $U$, but we do not investigate this.  
\end{remark}

\begin{proposition}
\label{proposition:uct_calculating_ei_cohomology_of_a_local_fp_type_collapses_at_a_finite_page}
Let $X$ be a locally fp spectrum and consider the universal coefficient spectral sequence
\begin{equation}
\label{equation:universal_coefficient_sseq_calculating_emodi_cohomology}
E_{2}^{\ast, \ast} = \Ext_{V_{*}}^{\ast}(V_{*}X, E_{*}/I) \Rightarrow (E \otimes M)^{*}(X). 
\end{equation}
Then (\ref{equation:universal_coefficient_sseq_calculating_emodi_cohomology}) collapses at a finite page, is strongly convergent, and $E^{\ast, \ast}_{r}$ is a finitely presented as an $\Ext_{V_{*}}^{\ast}(E_{*}/I, E_{*}/I)$-module for each $r \geq 2$. 
\end{proposition}

\begin{proof}
The spectral sequence of (\ref{equation:universal_coefficient_sseq_calculating_emodi_cohomology}) is a module over the spectral sequence 
\begin{equation}
\label{equation:base_ring_spectral_sequence}
\Ext_{V_{*}}^{\ast}(E_{*}/I, E_{*}/I) \Rightarrow \pi_{*} \map_{\Mod_{V}}(E \otimes M, E \otimes M),
\end{equation}
which collapses on the second page by \cref{lemma:collapse_of_uct_for_emodi_on_second_page}. By the assumption that $X$ is locally fp, $V_{*}X$ is degreewise finite and hence finitely presented as a module over $V_{*}$ by \cref{lemma:characterization_of_fp_modules_over_homotopy_of_generalized_moore}. Consequently,  \cref{lemma:ext_algebra_noether_and_finite_presentation_over_ext_algebra} implies that the $E_{2}$-page of (\ref{equation:universal_coefficient_sseq_calculating_emodi_cohomology}) is finitely presented as a module over the noetherian ring $\Ext_{V_{*}}^{\ast}(E_{*}/I, E_{*}/I)$. By a classical argument about a finitely generated module over a collapsing spectral sequence with a noetherian second page, see \cite[Proposition 4.1]{bajer1994may}, this implies that each further page is also finitely generated and that the spectral sequence collapses at a finite page. 

For convergence, as in the proof of \cref{lemma:collapse_of_uct_for_emodi_on_second_page} we observe that strong convergence follows from collapse at a finite page by a result of Boardman \cite[Remark following Theorem 7.1]{boardman1999conditionally}.
\end{proof} 

\begin{proof}[{Proof of \cref{theorem:classification_of_local_fp_types_in_terms_of_cohomology_and_morava_e_theory}, implication $(1 \Rightarrow 2)$:}]
Let $X$ be locally fp. By \cref{proposition:uct_calculating_ei_cohomology_of_a_local_fp_type_collapses_at_a_finite_page}, the universal coefficient spectral sequence of (\ref{equation:universal_coefficient_sseq_calculating_emodi_cohomology}) collapses at some finite page $E_{r}$, which is finitely generated over the $E_{2} = E_{r}$-page of (\ref{equation:base_ring_spectral_sequence}). In other words, 
if we look at the induced filtrations on the abutments (i.e. the $E_{\infty}$-term), then the associated graded of one is finitely generated as a module over the associated graded of the other. 

As both filtrations on abutments are complete and Hausdorff due to strong convergence, \cite[(I), Proposition 5.3]{nastasescu2006graded} shows that $(E \otimes M)^{\ast}(X)$ is finitely generated as a module over $(E^{\ast}/I) \llbracket U \rrbracket$, where we are using the identification of \cref{remark:self_maps_of_emodi_as_an_f_module}. In particular, $(E \otimes M)^{*}(X)$ is finitely generated over $E^{*} \llbracket U \rrbracket$ and hence also over
\[
E^{*} \llbracket \mathbf{G}_{n} \rrbracket \simeq \pi_{-*} \map_{\spectra}(E, E) \simeq E^{*}E.  
\]

If $M = \thesphere / (p^{\alpha_{0}}, v_{1}^{\alpha_{1}}, \ldots, v_{n-1}^{\alpha_{n-1}})$, we write 
\[
E/I_{k} \colonequals E \otimes \thesphere / (p^{\alpha_{0}}, \ldots, v_{k-1}^{\alpha_{k-1}}). 
\]
We will show by downward induction that 
\[
(E / I_{k})^{*}(X)
\]
is finitely generated over $E^{*}E$ for $0 \leq k \leq n$. The base  case of $k = n$ was covered in the previous paragraph. For the inductive step, observe that the long exact sequence of cohomology associated to multiplication by $v_{k}^{\alpha_{k}}$ shows that 
\[
((E/I_{k})^{*}X)/v_{k}^{\alpha_{k}} \hookrightarrow (E/I_{k+1})^{*}X
\]
is a monomorphism. By the inductive assumption, the target is finitely generated over $E^{*}E$, and hence so is the source, since $E^{*}E$ is noetherian by \cref{corollary:twisted_group_algebra_of_gn_is_noether}. As $(E/I_{k})^{*}X$ is derived $v_{k}$-complete as the homotopy of a $K(n)$-local $E$-module, the derived Nakayama lemma, see \cite[\href{https://stacks.math.columbia.edu/tag/09B9}{Tag 09B9}]{stacks-project}, implies that any set of generators of the quotient by $v_{k}^{\alpha_{k}}$ lifts to a set of generators of $(E/I_{k})^{*}(X)$. This completes the inductive step and shows finite generation as needed. 
\end{proof}

\section{The height one case of the Hahn-Wilson conjecture} 

In this short section, we use the characterization of locally fp spectra to deduce the original Hahn-Wilson conjecture at height one. This reduction is not particularly difficult, and relies on two crucial facts: 
\begin{enumerate}
    \item the telescope conjecture is true at height one, 
    \item if $X$ is of fp-type $1$, then $\pi_{\ast}(X) \rightarrow \pi_{\ast}(L_{K(1)} X)$ is an isomorphism in large degrees.  
\end{enumerate}

\begin{remark}
While the telescope conjecture is now known to fail in general at heights $n > 1$ \cite{burklund2023k}, the analogue of the second property holds for spectra of fp-type $n > 1$ if we replace $L_{K(1)}$-localization with $L_{n}^{f}$-localization \cite[{Theorem 3.1.3}]{hahn2022redshift}. 

This raises the interesting possibility of showing the Hahn-Wilson conjecture for spectra of fp-type n which satisfy the telescope conjecture by a suitable height decomposition and reduction to \cref{theorem:classification_of_local_fp_types_in_terms_of_cohomology_and_morava_e_theory}. We hope to return to this question in future work. 
\end{remark}

\begin{notation}
All spectra are implicitly $p$-complete, so that $\BPn{1}$ denotes the $p$-complete truncated Brown-Peterson spectrum. Since we specialize to height $1$, we use $E \colonequals E_{1}$ to denote the first Morava $E$-theory. In our setup, $E_1$ is the $p$-complete complex $K$-theory $\mathrm{KU}_p^\wedge$, and $\BP\langle1\rangle$ is the Adams summand $\ell_p^\wedge$ as a spectrum unambiguously \cite{adamsuniqueness,leeuniqueness}.
\end{notation}

\begin{lemma}
\label{lemma:thick_subcat_generated_by_1st_e_truncate_to_thick_subcat_of_bp1}
Let $X$ be any spectrum in the thick subcategory generated by $E$. Then, the connective cover $\tau_{\geq0}X$ is in the thick subcategory generated by $\BPn{1}$.
\end{lemma}

\begin{proof}
Let $\mathcal C$ denote the collection of all $p$-complete spectra $X$ such that $\pi_i(X)$ is finitely generated over $\mathbf Z_p$ for every $i$ and such that $\tau_{\geq 0} X$ is in the thick subcategory generated by $\BP\langle1\rangle$. We will show that 
\begin{enumerate}
    \item $E_1\in \mathcal C$, 
    \item $\mathcal C$ is a thick subcategory of the category of all $p$-complete spectra.
\end{enumerate}
Together, they imply the needed claim. The first statement follows from the fact that $\mathrm{KU}_p^\wedge$ splits into shifts of the Adams summand
\[
    \tau_{\geq0}E_1\simeq \BP\langle1\rangle\oplus\cdots\oplus\Sigma^{2p-4}\BP\langle1\rangle.
\]

For the second statement, observe that $\mathcal C$ is clearly closed under retracts. Suppose that $X\to Y\to Z$ is a cofiber sequence where $X,Y\in\mathcal C$. Then, there is a cofiber sequence
\[
    \cof(\tau_{\geq0}X\to\tau_{\geq0}Y)\to \tau_{\geq0}Z \to \Ker(\pi_{-1}(X)\to\pi_{-1}(Y))
\]
where the last term is considered as a discrete spectrum. Since $\mathbf Z_p$ is in the thick subcategory generated by $\BPn{1}$, any finitely generated $\mathbf Z_p$-module, considered as a discrete spectrum, is in the thick subcategory generated by $\BPn{1}$. Therefore, since $\Ker(\pi_{-1}(X)\to\pi_{-1}(Y))$ is a finitely generated $\mathbf Z_p$-module, $Z$ is contained in $\mathcal C$. An analogous argument shows closure under fibers.
\end{proof}

\begin{theorem}
\label{theorem:main_text_hahn_wilson_conjecture_at_height_one}
Let $X$ be a spectrum of fp-type $\leq1$. Then, $X$ is in the thick subcategory generated by $\BPn{1}$.
\end{theorem}
\begin{proof}
We may assume that $X$ is connective. Let $V$ be a type $1$ finite complex with $v_1$-self map $v \colon V \rightarrow V$. Since $V/v \otimes X$ is $\pi_{\ast}$-finite, in particular bounded, the fiber of
\[
    V \otimes X \to V[v^{-1}] \otimes X \simeq L_{K(1)}(V \otimes X)
\]
is bounded above. Here, the given equivalence follows from the truth of the telescope conjecture at height one  \cite{mahowald1981bo, miller1981relations}. We deduce that the fiber of $X\to L_{K(1)}X$ is bounded above. 

Since $L_{K(1)}X$ is in the thick subcategory generated by $E$ by \cref{theorem:classification_of_local_fp_types_in_terms_of_cohomology_and_morava_e_theory}, we have that $\tau_{\geq 0} (L_{K(1)}X)$ is in the thick subcategory generated by $\BPn{1}$ by \cref{lemma:thick_subcat_generated_by_1st_e_truncate_to_thick_subcat_of_bp1}. We deduce that the same is true for $X$ since the fiber of $X\to \tau_{\geq 0}(L_{K(1)}X)$ is bounded and has finitely generated homotopy groups.
\end{proof}


\appendix

\bibliographystyle{amsalpha}
\bibliography{fp_types_bibliography}

\end{document}